\documentclass[10pt,a4paper]{article}
\usepackage{amssymb}
\usepackage{latexsym}
\usepackage{lineno,hyperref}
\usepackage{hyperref}
\usepackage{appendix}

\hypersetup{colorlinks,linkcolor= blue,citecolor=blue}
\newtheorem{theorem}{Theorem}[section]

\newtheorem{lemma}{Lemma}[section]

\newtheorem{remark}{Remark}[section]

\def\proof{\mbox {\it \textbf{Proof.}~~}}
\usepackage{amssymb}
\usepackage{amsmath}
\numberwithin{equation}{section}
\begin{document}
\title{{\bf\Large  Variational methods for degenerate Kirchhoff equations}}
\author{{\small Sheng-Sen Lu}\\
{\it\small Center for Applied Mathematics}, {\it\small Tianjin University}\\
{\it\small Tianjin, 300072,}\\
{\it\small e-mail: sslu@tju.edu.cn}}
\date{}
\maketitle
\begin{center}
{\bf\small Abstract}\\
\vspace{2mm}
\hspace{.00in}\parbox{4.5in}
{{\small
For a degenerate autonomous Kirchhoff equation which is set on $\mathbb{R}^N$ and involves the Berestycki-Lions type nonlinearity, we cope with the cases $N=2,3$ and $N\geq5$ by using mountain pass and symmetric mountain pass approaches and by using Clark theorem respectively.}}
\end{center}
\noindent
{ \it \small Key words}: {\small Kirchhoff equations, Ground state, Multiple solutions, Nonexistence result,\\ Variational methods.}

\section{Introduction}\label{section 1}
In this paper, we continue the study initiated in \cite{Lu17} on the existence of a positive ground state (or least energy) solution and multiple solutions to the following autonomous Kirchhoff problem
\begin{equation*}\tag{$P$}\label{problem P}
\left\{
\begin{aligned}
&-\left(a+b\int_{\mathbb{R}^N}|\nabla{u}|^2\right)\Delta{u}= f(u)~~\text{in}~\mathbb{R}^N,\\
&u\in H^1(\mathbb{R}^N),~~~~u\not\equiv0~~\text{in}~\mathbb{R}^N,
\end{aligned}
\right.
\end{equation*}
where $a\geq0,b>0$ are constants, $N\geq1$ and $f:\mathbb{R}\rightarrow\mathbb{R}$ is a given function that satisfies the so-called Berestycki-Lions type conditions.

Such problems and their variations have been extensively studied especially in the recent decade. The interests are mainly originated from two aspects: one is the strong physical background they come from, and the other is the presence of some specific mathematical difficulties that make their study challenging. For more details, we refer readers to \cite{Az11,Az12,Az15,Fi14,Lu15,Lu17,Lu16} and the references therein.

Now we recall the Berestycki-Lions type conditions, which are assumed on the nonlinearity $f$ and specified as follows:
\begin{itemize}
  \item[$(f1)$] $f\in C(\mathbb{R},\mathbb{R})$ is continuous and odd.
  \item[$(f2)$] $-\infty<\underset{t\rightarrow0}{\liminf}f(t)/t\leq\underset{t\rightarrow0}{\limsup}f(t)/t<0$.
  \item[$(f3)$] When $N\geq3$, $\underset{t\rightarrow\infty}{\lim}f(t)/|t|^{\frac{N+2}{N-2}}=0$.

  When $N=2$,
  $\underset{t\rightarrow\infty}{\lim}f(t)/e^{\alpha t^2}=0$ for any given $\alpha>0$,
  \item[$(f4)$] Let $F(t):=\int^t_0f(\tau)d\tau$. When $N\geq2$, there exists $\zeta>0$ such that $F(\zeta)>0$.

  When $N=1$, there exists $\zeta>0$ such that
  \begin{equation*}
  F(t)<0~\text{for all}~t\in (0,\zeta),~~F(\zeta)=0~~\text{and}~~f(\zeta)>0.
  \end{equation*}
\end{itemize}
Such type of conditions were first introduced by Berestycki and Lions \cite{Be83-1,Be83-2} in the study of the following nonlinear scalar field equation
\begin{equation}\tag{$Q$}\label{problem Q}
\left\{
\begin{aligned}
&-\Delta{v}= f(v)~~&\text{in}~\mathbb{R}^N,\\
&v\in H^1(\mathbb{R}^N),~~v\not\equiv0~~&\text{in}~\mathbb{R}^N,
\end{aligned}
\right.
\end{equation}
After that, some generalizations were made in \cite{Be83-3,Je03,Hi10,Fi14}. We remark that the conditions $(f1)-(f4)$ are almost necessary and turn out to be sufficient to get a nontrivial solution of Problem \eqref{problem Q}. When it comes to Problem \eqref{problem P}, these conditions are still almost necessary but may not be sufficient, see e.g. \autoref{thm 1.3} below.

Under the very general conditions on $f$ as above, we investigated Problem \eqref{problem P} in \cite{Lu17} for all the cases $a\geq0$, $b>0$ and $N\geq1$. In that paper, by using certain known results about the solutions to Problem \eqref{problem Q} and a scaling argument due to Azzollini \cite{Az12,Az15}, a positive ground state solution and multiple solutions to Problem \eqref{problem P} were obtained for $N\geq2$. When $N=1$, the existence and the uniqueness (in a certain sense) of the nontrivial solutions were shown. In addition, certain nonexistence results and other interesting ones were also observed. For further details, we refer readers to \cite{Lu17}. See also \cite{Az11,Az12,Az15,Fi14,Lu16}, which consider Problem \eqref{problem P} in the non-degenerate case $a>0$.

As one can see, the main proofs in \cite{Lu17} are non-variational. One of the key points there is the observation that solutions of the nonlocal Problem \eqref{problem P} can be characterized as that of the local Problem \eqref{problem Q} with a suitable spatial scaling. On the other hand, in view of $(f1)-(f3)$, solutions of Problem \eqref{problem P} can also be naturally characterized as critical points of the functional $\Phi\in C^1(H^1(\mathbb{R}^N),\mathbb{R})$ defined by
\begin{equation*}
\Phi(u)=\frac{1}{2}a\int_{\mathbb{R}^N}|\nabla{u}|^2+\frac{1}{4}b\left(\int_{\mathbb{R}^N}|\nabla{u}|^2\right)^2-\int_{\mathbb{R}^N}F(u).
\end{equation*}
Thus, it would be interesting to know whether for this general Berestycki-Lions type nonlinearity $f$ one can still obtain the main results of \cite{Lu17} via variational methods.

Up to now, this question has already been answered positively and almost completely in the non-degenerate case $a>0$, see the early papers \cite{Az12,Az11,Fi14} and the more recent one \cite{Lu16}. While, to the best of our knowledge, the degenerate case of this question is totally open and seems to be the very challenging part.

Motivated by the open question raised above and enlightened by the nice results in \cite{Lu17}, in the present paper, we are going to mainly deal with the degenerate case of Problem \eqref{problem P} via variational methods. To be more precise, in the case $a=0, b>0$, we will show the existence of infinitely many distinct radial solutions for $N=2,3$ and $N\geq5$ via symmetric mountain pass argument and Clark theorem respectively. With the aid of mountain pass approach, a positive ground state solution will also be obtained in the case $a=0,b>0$ and $N=2,3$. In addition, we will also treat the non-degenerate case $a,b>0,N\geq5$ by using Clark theorem and show the existence of multiple radial solutions with negative energies for suitable $a,b>0$.

Now, we state our main results of this paper as follows.

\begin{theorem}\label{thm 1.1}
Assume that $a=0,b>0$ fixed, $N=2,3$ and $f$ satisfies $(f1)-(f4)$. Then Problem \eqref{problem P} has a positive ground state (or least energy) solution $u_0$ and infinitely many distinct radial solutions $\{u_k\}^{+\infty}_{k=1}$ for any $b>0$. Moreover, $u_0$ is radially symmetric, $\Phi(u_0)>0$ and $\Phi(u_k)\to+\infty$ as $k\to+\infty$.
\end{theorem}
\begin{remark}\label{remark 1.1}
\begin{description}
  \item[$~~(i)$] \autoref{thm 1.1} was first established in \cite{Lu17} by using a non-variational method, precisely, a scaling argument. As shown in Section \ref{section 3}, we can also prove this result via mountain pass and symmetric mountain pass approaches. Thus, for this degenerate case, we answer the open question raised above in the affirmative.
  \item[$~(ii)$] It is worth pointing out that, in this case, we also provide in this paper a mountain pass characterization of ground state solutions of Problem \eqref{problem P}, see \autoref{lemma 3.5} and \autoref{thm C.1}. Such a characterization is expected to be useful in the studies of the corresponding non-autonomous and singular perturbation problems.
\end{description}
\end{remark}

\begin{theorem}\label{thm 1.2}
Assume that $a\geq0$ fixed, $b>0$, $N\geq5$ and $f$ satisfies $(f1)-(f4)$. Then the following statements hold.
\begin{description}
  \item[$~~(i)$] If $a>0$, then for any given positive integer $k$ there exists a positive constant $b_k>0$ such that Problem \eqref{problem P} has at least $k$ distinct radial solutions with negative energies for any $b\in(0,b_k)$.
  \item[$~(ii)$] If $a=0$, then Problem \eqref{problem P} has infinitely many distinct radial solutions $\{u_k\}^{+\infty}_{k=1}$ for any $b>0$. Moreover, $\Phi(u_k)<0$ for all $k$ and $u_k\to0$ in $H^1(\mathbb{R}^N)$ as $k\to+\infty$.
\end{description}
\end{theorem}
\begin{remark}\label{remark 1.2}
\begin{description}
  \item[$~~(i)$] In the case $a>0$ fixed and $N\geq5$, multiple radial solutions with positive energies were first shown in \cite{Az11} (for sufficiently small $b>0$) via symmetric mountain pass approach and a truncation argument. See also the recent work \cite{Lu16} for a slight different proof, which does not involve the truncation technique. Then, with the aid of a scaling argument, a progress was made in \cite{Lu17} which claims the existence of multiple radial solutions with negative energies. In the present paper, we manage to give a variational proof of this progress by using Clark theorem, see Item $(i)$ of \autoref{thm 1.2} above and its proof in Section \ref{section 4}.
  \item[$~(ii)$]  In the case $a=0,b>0$ fixed and $N\geq5$, we know from Poho\u{z}aev identity that Problem \eqref{problem P} has no nontrivial solutions with nonnegative energies. Thus, it is natural to try to find solutions with negative energies. Based on a non-variational argument, it was proved in \cite{Lu17} that Problem \eqref{problem P} has infinitely many distinct radial solutions the energies of which are all negative and converge to $0$. In the present paper, by using Clark theorem, we not only obtain the same multiplicity result but also show the new fact that $\{u_k\}^{+\infty}_{k=1}$ converges strongly to zero in $H^1(\mathbb{R}^N)$, see Item $(ii)$ of \autoref{thm 1.2} above and its proof in Section \ref{section 4}. For another interesting fact, we refer readers to \autoref{remark 4.2}.
  \item[$(iii)$] When $a\geq0,b>0$ and $N\geq5$, some variational arguments in \cite{Az15,Lu17} showed that the global infimum $m_{\inf}$ of $\Phi$ is finite and can be achieved by a nonnegative radial function (but being not necessarily nonzero). Clearly, at least for the case $m_{\inf}<0$, we can say that a positive ground state solution has actually been obtained in \cite{Az15,Lu17} (via variational methods). It is already the case for sufficiently small $b>0$ if $a>0$ fixed and for all $b>0$ if $a=0$, see e.g. \autoref{lemma 2.2} below. This result is essentially due to \cite{Az15,Lu17} and we highlight it here just for completeness.
  \item[$(iv)$] The reader should be aware that, in the case $a>0$ fixed and $N\geq5$, the global infimum $m_{\inf}$ always be zero if $b$ is not less than a certain constant $b^{**}>0$, see Item $(iii)$ of Theorem 1.2 in \cite{Lu17}. But, there exists a positive ground state solution with positive or zero energy if $b$ is also not more than another certain constant $b^*>b^{**}$, see Remark 3.3 in \cite{Lu17}. It is very interesting to know how to cope with this subtle case by variational methods. By now, this question still remains open.
\end{description}
\end{remark}

In addition to the above existence and multiplicity results, we also prove in this paper the following nonexistence result.
\begin{theorem}\label{thm 1.3}
Assume that $f$ satisfies $(f1)-(f4)$, $N\geq4$ if $a>0$ fixed, $N=4$ if $a=0$. Then there exists a positive constant $b_*>0$ (dependent only on $a,N$ and $f$) such that Problem \eqref{problem P} has no nontrivial solutions when $b>b_*$.
\end{theorem}
\begin{remark}\label{remark 1.3}
\autoref{thm 1.3} is actually a special case of the sharp nonexistence results proved in \cite{Az15,Lu17}. But, the proof we provide here is variational and still works in the non-autonomous case (under suitable assumptions). Since the proof is short, we give it here. When $N\geq4$, in view of $(f1)-(f3)$, a positive constant $C_f>0$ exists such that $F(t)\leq C_f |t|^{2N/(N-2)}$ (or $f(t)t \leq C_f |t|^{2N/(N-2)}$) for all $t\in\mathbb{R}$. Thus, for any given solution $u$ of Problem \eqref{problem P}, we know from Poho\u{z}aev identity (or the fact that $\Phi'(u)u=0$) and Sobolev inequality that
\begin{equation*}
a\int_{\mathbb{R}^N}|\nabla{u}|^2+b\left({\int_{\mathbb{R}^N}}|\nabla{u}|^2\right)^2\leq C_f C_N\left({\int_{\mathbb{R}^N}}|\nabla{u}|^2\right)^{\frac{N}{N-2}}.
\end{equation*}
Assume $N\geq4$ if $a>0$ fixed, $N=4$ if $a=0$, then we can find a sufficiently large positive constant $b_*>0$ such that
\begin{equation*}
at+b_*t^2>C_f C_Nt^{\frac{N}{N-2}}~~~~\text{for all}~t>0.
\end{equation*}
Obviously, the solution $u$ must be zero if $b>b_*$.
\end{remark}

Now, we can say that the open question raised above has been answered positively and satisfactorily in most degenerate cases. More important, the present paper provides us with a better understanding about these degenerate cases from a variational point of view. This would be helpful to the study of the corresponding non-autonomous case (by using variational methods), where the scaling argument used in \cite{Lu17} does not work in general.

The remaining part of this paper is organized as follows. In Section \ref{section 2}, we show some geometrical properties of the functional $\Phi$ and an important compactness result. With the aid of minimax methods, \hyperref[thm 1.1]{Theorems \ref{thm 1.1}} and \ref{thm 1.2} are proved in Sections \hyperref[section 3]{\ref{section 3}} and \ref{section 4} respectively. Appendix \hyperref[appendix]{A} involves some more technical auxiliary results that are crucial to conduct the proof of \autoref{thm 1.1} in Section \ref{section 3}.

\section{Preliminary}\label{section 2}
In this section, we shall establish some useful preliminary results for Problem \eqref{problem P} which will be exploited later in the variational proofs of the main results.

\subsection{Some geometrical properties of the functional}\label{subsection 2.1}

In what follows, for every positive integer $k\in\mathbb{N}^+$, we set
$$\mathbb{D}_k:=\left\{\sigma=(\sigma_1,\cdots,\sigma_k)\in \mathbb{R}^k~|~|\sigma|\leq1\right\}~~~~\text{and}~~~~\mathbb{S}^{k-1}:=\partial\mathbb{D}_k.$$
Also, for convenience, let
\begin{equation*}
\nu:=-\frac{1}{2}\underset{t\rightarrow0}{\limsup}\frac{f(t)}{t}\in (0,+\infty).
\end{equation*}

When $N=2,3$, it has been shown in \cite{Lu16} that the natural functional $\Phi$ corresponding to Problem \eqref{problem P} has the symmetric mountain pass geometry in the non-degenerate case $a,b>0$. Actually, as we can see below, the functional $\Phi$ still has such a geometry even in the degenerate case we consider in this paper. This fact seems not to be clearly known in the literature.
\begin{lemma}\label{lemma 2.1}
Assume that $a=0$, $b>0$ fixed, $N=2,3$ and $f$ satisfies $(f1)-(f4)$. Then the functional $\Phi$ satisfies the following properties.
\begin{itemize}
  \item[~~$(i)$] There exist $r_0>0$ and $\rho_0>0$ such that
                 \begin{equation*}
                 \begin{split}
                        \Phi(u)>0~~~~&\text{for any}~u\in H^1(\mathbb{R}^N)~\text{with}~0<\|u\|_{H^1}\leq r_0,\\
                        \Phi(u)\geq\rho_0~~~&\text{for any}~u\in H^1(\mathbb{R}^N)~\text{with}~\|u\|_{H^1}= r_0.
                 \end{split}
                 \end{equation*}
  \item[~$(ii)$] For every $k\in \mathbb{N}^+$, there exists an odd continuous mapping $\gamma_{0k}:\mathbb{S}^{k-1}\to H^1_r(\mathbb{R}^N)\setminus\{0\}$ such that
    \begin{equation*}
    \max_{\sigma\in \mathbb{S}^{k-1}}\Phi(\gamma_{0k}(\sigma))<0.
    \end{equation*}
\end{itemize}
\end{lemma}
\proof $(i)$ To prove this item, we mainly adopt the argument explored in the proof of Lemma 1.1 of \cite{Je02} but with certain technical modifications since we are now dealing with the degenerate Kirchhoff problem. First, we prove the case $N=3$. For $\nu>0$ defined above, from $(f2)$ and $(f3)$, a positive constant $C_{\nu}>0$ can be found such that
\begin{equation*}
-f(t)\geq \nu t-C_{\nu}t^5~~~~\text{for all}~t\geq0.
\end{equation*}
Then, by $(f1)$, we have
\begin{equation*}
-F(t)\geq \frac{1}{2}\nu t^2-\frac{1}{6}C_{\nu}t^6~~~~\text{for all}~t\in\mathbb{R}.
\end{equation*}
In view of the embedding $H^1(\mathbb{R}^3)\hookrightarrow L^6(\mathbb{R}^3)$, a positive constant $C'_\nu>0 $ exists such that, for any $u\in H^1(\mathbb{R}^3)$ with $||u||_{H^1}\leq1$,
\begin{equation*}
\begin{aligned}
\Phi(u)&\geq \frac{1}{4}b\left(\int_{\mathbb{R}^3}|\nabla u|^2\right)^2+\frac{1}{2}\nu\int_{\mathbb{R}^3}u^2-\frac{1}{6}C_\nu\int_{\mathbb{R}^3}u^6\\
&\geq\frac{1}{8}\min\{b,\nu\}\left(\int_{\mathbb{R}^3}|\nabla u|^2+\int_{\mathbb{R}^3}u^2\right)^2-\frac{1}{6}C_\nu\int_{\mathbb{R}^3}u^6\\
&\geq \frac{1}{8}\min\{b,\nu\}\|u\|^4_{H^1}-\frac{1}{6}C'_\nu\|u\|^6_{H^1}.
\end{aligned}
\end{equation*}
Now we can complete the proof of this case by choosing $r_0,\rho_0\in(0,1)$ sufficiently small .

Next, we deal with the remaining case $N=2$. For any given $\alpha>0$, by $(f2)$ and $(f3)$, a positive constant $C_\alpha>0$ exists such that
\begin{equation*}
-f(t)\geq \nu t-C_\alpha t^7 e^{\alpha t^2}~~~~\text{for all}~t\geq0.
\end{equation*}
Since
\begin{equation*}
\int^t_0\tau^7e^{\alpha \tau^2} d\tau=\frac{1}{2\alpha}t^6\left(e^{\alpha t^2}-1\right)-\frac{3}{\alpha}\int^t_0\tau^5\left(e^{\alpha \tau^2}-1\right)d\tau\leq \frac{1}{2\alpha}t^6\left(e^{\alpha t^2}-1\right)
\end{equation*}
for all $t\geq0$ and $f$ is an odd continuous function, we have
\begin{equation*}
-F(t)\geq \frac{1}{2}\nu t^2-\frac{C_\alpha}{2\alpha}t^6\left(e^{\alpha t^2}-1\right)~~~~\text{for all}~t\in\mathbb{R}.
\end{equation*}
Then, from the embedding $H^1(\mathbb{R}^2)\hookrightarrow L^{12}(\mathbb{R}^2)$, a positive constant $C'_\alpha>0$ can be found such that, for any $u\in H^1(\mathbb{R}^2)$ with $\|u\|_{H^1}\leq1$,
\begin{equation*}
\begin{aligned}
\Phi(u)&\geq \frac{1}{4}b\left(\int_{\mathbb{R}^2}|\nabla u|^2\right)^2+\frac{1}{2}\nu\int_{\mathbb{R}^2}u^2-\frac{C_\alpha}{2\alpha}\int_{\mathbb{R}^2}u^6\left(e^{\alpha u^2}-1\right)\\
&\geq \frac{1}{8}\min\{b,\nu\}\left(\int_{\mathbb{R}^2}|\nabla u|^2+\int_{\mathbb{R}^2}u^2\right)^2-\frac{C_\alpha}{2\alpha}\left(\int_{\mathbb{R}^2}u^{12}\right)^{\frac{1}{2}}\left[\int_{\mathbb{R}^2}\left(e^{\alpha u^2}-1\right)^2\right]^{\frac{1}{2}}\\
&\geq \frac{1}{8}\min\{b,\nu\}\|u\|^4_{H^1}-\frac{C'_\alpha}{2\alpha}\|u\|^6_{H^1}\left[\int_{\mathbb{R}^2}\left(e^{2\alpha u^2}-1\right)\right]^{\frac{1}{2}}.
\end{aligned}
\end{equation*}
We also know from the Moser-Trudinger inequality (see e.g. \cite{Ad00,Og90}) that there exist $\alpha_0>0$ and $M>0$ such that
\begin{equation*}
\int_{\mathbb{R}^2}\left(e^{\alpha_0 u^2}-1\right)\leq M~~~~\text{for all}~\|u\|_{H^1}\leq1.
\end{equation*}
Thus, by setting $\alpha:=\frac{1}{2}\alpha_0>0$ and choosing $r_0,\rho_0\in(0,1)$ sufficiently small, we can now complete the proof of this remaining case.

\smallskip
$(ii)$ For any $N\geq2$ and every $k\in \mathbb{N}^+$, arguing as in Theorem 10 of \cite{Be83-2}, an odd continuous mapping $\pi_k:\mathbb{S}^{k-1}\to H^1_r(\mathbb{R}^N)$ is defined such that
\begin{equation*}
0\notin \pi_k(\mathbb{S}^{k-1})~~~~\text{and}~~~~\underset{\sigma\in\mathbb{S}^{k-1}}{\min}\int_{\mathbb{R}^N}F(\pi_k(\sigma))\geq1.
\end{equation*}
Then a positive constant $\alpha_k>0$ exists such that
\begin{equation*}
\underset{\sigma\in\mathbb{S}^{k-1}}{\max}\int_{\mathbb{R}^N}|\nabla \pi_k(\sigma)|^2\leq \alpha_k.
\end{equation*}
For any $\sigma\in\mathbb{S}^{k-1}$, setting $\beta^s_k(\sigma)(x):=\pi_k(\sigma)(s^{-1}x)$ with $s>0$ undetermined, we have
\begin{equation*}
\begin{split}
\Phi(\beta^s_k(\sigma))
&=\frac{1}{2}as^{N-2}\int_{\mathbb{R}^N}|\nabla{\pi_k(\sigma)}|^2\\
&~~~~~~~~~~~~~~+\frac{1}{4}bs^{2N-4}\left(\int_{\mathbb{R}^N}|\nabla{\pi_k(\sigma)}|^2\right)^2-s^N\int_{\mathbb{R}^N}F(\pi_k(\sigma))\\
&\leq\frac{1}{2}a\alpha_ks^{N-2}+\frac{1}{4}b\alpha^2_ks^{2N-4}-s^N=:g_k(s).
\end{split}
\end{equation*}
When $a=0$, $b>0$ and $N=2,3$, it is clear that $g_k(s_k)<0$ for sufficiently large $s_k>0$. Thus, the proof of this item is completed by defining $\gamma_{0k}:=\beta^{s_k}_k$.~~$\square$

Concerning the case where $a\geq0,b>0$ and $N\geq5$, we have the following result which is also related to the geometrical properties of $\Phi$ and is necessary when we try to deal with this case via Clark theorem.
\begin{lemma}\label{lemma 2.2}
Assume that $a\geq0$, $b>0$ fixed, $N\geq5$ and $f$ satisfies $(f1)-(f4)$. Then the following statements hold.
\begin{description}
   \item[$~~(i)$] The functional $\Phi$ is bounded from below and coercive with respect to $H^1\text{-norm}$.
   \item[$~(ii)$] If $a>0$, then for every $k\in\mathbb{N}^+$ there exist a positive constant $b_k>0$ and an odd continuous mapping $\overline{\gamma}_{0k}:\mathbb{S}^{k-1}\to H^1_r(\mathbb{R}^N)\setminus\{0\}$ such that, for $b\in(0,b_k)$,
       \begin{equation*}
    \underset{\sigma\in \mathbb{S}^{k-1}}{\max}\Phi(\overline{\gamma}_{0k}(\sigma))<0.
    \end{equation*}
   \item[$(iii)$] If $a=0$, then for any fixed $b>0$ and every $k\in \mathbb{N}^+$ there exists an odd continuous mapping $\widetilde{\gamma}_{0k}:\mathbb{S}^{k-1}\to H^1_r(\mathbb{R}^N)\setminus\{0\}$ such that
    \begin{equation*}
    \underset{\sigma\in \mathbb{S}^{k-1}}{\max}\Phi(\widetilde{\gamma}_{0k}(\sigma))<0.
    \end{equation*}
\end{description}
\end{lemma}
\proof $(i)$ This item has actually been proved in \cite{Az15} and \cite{Lu17} for the cases $a>0$ and $a=0$ respectively. But, for reader's convenience, we provide a detailed proof here. When $N\geq5$, by $(f2)$ and $(f3)$, a positive constant $C_{\nu}>0$ exists such that
\begin{equation*}
-f(t)\geq \nu t-C_{\nu}t^{\frac{N+2}{N-2}}~~~~\text{for all}~t\geq0.
\end{equation*}
Then, since $f$ is an odd continuous function, we have
\begin{equation*}
-F(t)\geq \frac{1}{2}\nu t^2-\frac{N-2}{2N}C_{\nu}|t|^{\frac{2N}{N-2}}~~~~\text{for all}~t\in\mathbb{R}.
\end{equation*}
In view of the embedding $\mathcal{D}^{1,2}(\mathbb{R}^N)\hookrightarrow L^{\frac{2N}{N-2}}(\mathbb{R}^N)$, a positive constant $C'_\nu>0 $ can be found such that, for any $u\in H^1(\mathbb{R}^N)$,
\begin{equation*}
\begin{aligned}
\Phi(u)&\geq \frac{1}{4}b\left(\int_{\mathbb{R}^N}|\nabla u|^2\right)^2+\frac{1}{2}\nu\int_{\mathbb{R}^N}u^2-\frac{N-2}{2N}C_\nu\int_{\mathbb{R}^N}|u|^{\frac{2N}{N-2}}\\
&\geq \frac{1}{4}b\left(\int_{\mathbb{R}^N}|\nabla u|^2\right)^2+\frac{1}{2}\nu\int_{\mathbb{R}^N}u^2-\frac{N-2}{2N}C'_\nu\left(\int_{\mathbb{R}^N}|\nabla u|^2\right)^{\frac{N}{N-2}}.
\end{aligned}
\end{equation*}
Observing that $1<\frac{N}{N-2}<2$ for $N\geq5$, we can now conclude the desired result.

To prove Items $(ii)$ and $(iii)$, we shall use again the odd continuous mapping $\beta^s_k$, the positive constant $\alpha_k>0$ and the function $g_k$, which are introduced in the proof of \autoref{lemma 2.1} and are actually still valid here.

\smallskip
$(ii)$ When $a>0$ fixed and $N\geq5$, for every $k\in\mathbb{N}^+$, let
\begin{equation*}
\overline{s}_k:=\sqrt{2a\alpha_k}>0~~~~\text{and}~~~~b_k:=\alpha^{-2}_k(\overline{s}_k)^{4-N}>0.
\end{equation*}
It is clear that $g_k(\overline{s}_k)<0$ for any fixed $b\in(0,b_k)$. Thus, the positive constant $b_k>0$ defined above and the odd continuous mapping $\overline{\gamma}_{0k}:=\beta^{\overline{s}_k}_k$ are the desired ones.

\smallskip
$(iii)$ When $a=0$, $b>0$ fixed and $N\geq5$, it is clear that $g_k(\widetilde{s}_k)<0$ for sufficiently small $\widetilde{s}_k>0$. Thus, the proof of this item is completed by defining $\widetilde{\gamma}_{0k}:=\beta^{\widetilde{s}_k}_k$~~$\square$

\subsection{A compactness result}\label{subsection 2.2}

It is clear that $H^1_r(\mathbb{R}^N):=\{u\in H^1(\mathbb{R}^N)|u(x)=u(|x|)\}$ is a natural constraint to look for critical points of $\Phi$, namely critical points of the functional $\Phi$ restricted to $H^1_r(\mathbb{R}^N)$ are true critical points in $H^1(\mathbb{R}^N)$. On the other hand, as we can see below, we can recover some compactness by choosing $H^1_r(\mathbb{R}^N)$ as the working space when $N\geq2$. Thus, from now on, we will directly define $\Phi$ on $H^1_r(\mathbb{R}^N)$ unless it is explicitly stated otherwise. In analogy with the well-known compactness result in \cite{Be83-2}, we have the following result.
\begin{lemma}\label{lemma 2.5}
Assume that $a\geq0$, $b>0$, $N\geq2$ and $f$ satisfies $(f1)-(f3)$. Then each bounded Palais-Smale sequence $\{u_n\}^{+\infty}_{n=1}\subset H^1_r(\mathbb{R}^N)$ of $\Phi$ has a convergent subsequence.
\end{lemma}

To prove \autoref{lemma 2.5}, we need the following two lemmas.
\begin{lemma}\label{lemma 2.3}
(\cite{Be83-1,St77}) For any $N\geq2$, the following conclusions hold.
\begin{description}
  \item[$~(i)$] $H^1_r(\mathbb{R}^N)\subset C(\mathbb{R^N}\setminus\{0\})$ and there exists a positive constant $C_N>0$ such that
                     \begin{equation*}
                       |u(x)|\leq C_N |x|^{-\frac{N-1}{2}}\|u\|_{H^1_r(\mathbb{R}^N)}
                     \end{equation*}
                for $u\in H^1_r(\mathbb{R}^N)$ and $|x|\geq1$.
  \item[$(ii)$] The embedding $H^1_r(\mathbb{R}^N)\hookrightarrow L^q(\mathbb{R}^N)$ is compact for $q\in(2,2^*)$, where $2^*:=\frac{2N}{N-2}$ if $N\geq3$ and $2^*:=+\infty$ if $N=2$.
\end{description}
\end{lemma}
\begin{lemma}\label{lemma 2.4}
(\cite{Az11,Be83-1,St77}) Let $P$ and $Q:\mathbb{R}\to\mathbb{R}$ be two continuous functions satisfying
\begin{equation*}
\underset{|t|\to+\infty}{\lim}\frac{P(t)}{Q(t)}=0,
\end{equation*}
$\{v_n\}^{+\infty}_{n=1},v$ and $z$ be measurable functions from $\mathbb{R}^N$ to $\mathbb{R}$ such that
\begin{equation*}
\underset{n\in\mathbb{N}}{\sup}\int_{\mathbb{R}^N}|Q(v_n(x))z(x)|dx<+\infty
\end{equation*}
and
\begin{equation*}
P(v_n(x))\to v(x)~~\text{a.e. in}~\mathbb{R}^N,~~~~\text{as}~n\to+\infty.
\end{equation*}
Then, for any bounded Borel set $\Omega\subset\mathbb{R}^N$,
\begin{equation*}
\int_{\Omega}|\left(P(v_n)-v(x)\right)z(x)|dx\to0~~~~\text{as}~n\to+\infty.
\end{equation*}
Moreover, if there also hold
\begin{equation*}
\underset{t\to0}{\lim}\frac{P(t)}{Q(t)}=0~~~~\text{and}~~~~\underset{|x|\to+\infty}{\lim}\left(\underset{n\in\mathbb{N}}{\sup}~|v_n(x)|\right)=0,
\end{equation*}
then we can conclude further that
\begin{equation*}
\int_{\mathbb{R}^N}|(P(v_n)-v(x))z(x)|dx\to0~~~~\text{as}~n\to+\infty.
\end{equation*}
\end{lemma}

\medskip
\noindent
\textbf{Proof of \autoref{lemma 2.5}.}~~Let $\{u_n\}^{+\infty}_{n=1}\subset H^1_r(\mathbb{R}^N)$ be a bounded Palais-Smale sequence of $\Phi$. Then, from Item $(ii)$ of \autoref{lemma 2.3}, we may assume that there exist a radial function $u\in H^1_r(\mathbb{R}^N)$ and a nonnegative constant $A\geq0$ such that, up to a subsequence, as $n\to+\infty$,
\begin{equation*}
\left\{
\begin{aligned}
  &~u_n\rightharpoonup u&&\text{in}~H^1_r(\mathbb{R}^N),\\
  &~u_n\to u & &\text{in}~L^q(\mathbb{R}^N),~\forall~q\in(2,2^*),\\
  &~u_n(x)\to u(x)&&\text{a.e. in}~\mathbb{R}^N,\\
  &\int_{\mathbb{R}^N}|\nabla u_n|^2\to A^2&&\text{in}~\mathbb{R}.
\end{aligned}
\right.
\end{equation*}

Set $P(t):=f(t)$, $Q(t):=|t|^{\frac{N+2}{N-2}}$ if $N\geq3$ and $Q(t):=e^{\alpha t^2}-1$ if $N=2$, $v_n:=u_n$, $v:=f(u)$ and $z\in C^\infty_0(\mathbb{R}^N)\cap H^1_r(\mathbb{R}^N)$. Here $\alpha>0$ is chosen so that
\begin{equation*}
\underset{n\in\mathbb{N}}{\sup}\int_{\mathbb{R}^2} |Q(u_n(x))|^2dx<+\infty,
\end{equation*}
which implies
\begin{equation*}
\underset{n\in\mathbb{N}}{\sup}\int_{\mathbb{R}^2}|Q(u_n(x))z(x)|dx<+\infty.
\end{equation*}
This is possible from the boundedness of $\{\|u_n\|_{H^1_r}\}^{+\infty}_{n=1}$ and the Moser-Trudinger inequality (see e.g. \cite{Ad00,Og90}). Thus, in view of $(f2)$, $(f3)$ and the boundedness of $\{\|u_n\|_{H^1_r}\}^{+\infty}_{n=1}$, \autoref{lemma 2.4} shows that
\begin{equation}\label{equ2.1}
\underset{n\to+\infty}{\lim}\int_{\mathbb{R}^N}f(u_n)z=\int_{\mathbb{R}^N}f(u)z.
\end{equation}
As a direct consequence, $u\in H^1_r(\mathbb{R}^N)$ is a weak solution of the following problem
\begin{equation*}
-\left(a+bA^2\right)\Delta{v}= f(v),~~~~v\in H^1_r(\mathbb{R}^N),
\end{equation*}
and then the following equality holds
\begin{equation}\label{equ2.2}
\left(a+bA^2\right)\int_{\mathbb{R}^N}|\nabla u|^2= \int_{\mathbb{R}^N}f(u)u.
\end{equation}

On the other hand, let
\begin{equation*}
f_1(t):=\left\{
\begin{aligned}
&\max\{f(t)+2\nu t,0\}&\text{for}~t\geq0,\\
&\min\{f(t)+2\nu t,0\}&\text{for}~t<0,
\end{aligned}
\right.
\end{equation*}
and
\begin{equation*}
f_2(t):=f_1(t)-f(t)~~~~\text{for}~t\in\mathbb{R},
\end{equation*}
where $\nu>0$ is a positive constant defined at the beginning of this section. Obviously, $f_1$ and $f_2$ are odd continuous function satisfying $(f3)$, $\lim_{t\rightarrow0}f_1(t)/t=0$, $f_2(t)t\geq 2\nu t^2$ for all $t\in\mathbb{R}$ and $0<2\nu=\liminf_{t\rightarrow0}f_2(t)/t\leq\limsup_{t\rightarrow0}f_2(t)/t<+\infty$. Set $P(t):=f_1(t)t$, $Q(t):=t^2+|t|^{\frac{2N}{N-2}}$ if $N\geq3$ and $Q(t):=e^{\alpha t^2}-1$ if $N=2$, $v_n:=u_n$, $v:=f_1(u)u$ and $z:=1$. Here $\alpha>0$ is chosen so that
\begin{equation*}
\underset{n\in\mathbb{N}}{\sup}\int_{\mathbb{R}^2} |Q(u_n(x))|dx<+\infty,
\end{equation*}
and this is possible from the boundedness of $\{\|u_n\|_{H^1_r}\}^{+\infty}_{n=1}$ and the Moser-Trudinger inequality (see e.g. \cite{Ad00,Og90}). Then, from Item $(i)$ of \autoref{lemma 2.3}, the boundedness of $\{\|u_n\|_{H^1_r}\}^{+\infty}_{n=1}$ and \autoref{lemma 2.4}, we know that
\begin{equation}\label{equ2.3}
\underset{n\to+\infty}{\lim}\int_{\mathbb{R}^N}f_1(u_n)u_n=\int_{\mathbb{R}^N}f_1(u)u.
\end{equation}
Besides, by Fatou's lemma, we have
\begin{equation*}
\underset{n\to+\infty}{\liminf}\int_{\mathbb{R}^N}f_2(u_n)u_n\geq\int_{\mathbb{R}^N}f_2(u)u.
\end{equation*}
Thus,
\begin{equation*}
\begin{aligned}
\int_{\mathbb{R}^N}f(u)u&=\int_{\mathbb{R}^N}f_1(u)u-\int_{\mathbb{R}^N}f_2(u)u\geq\underset{n\to+\infty}{\lim}\int_{\mathbb{R}^N}f_1(u_n)u_n-\underset{n\to+\infty}{\liminf}\int_{\mathbb{R}^N}f_2(u_n)u_n\\
&\geq\underset{n\to+\infty}{\lim}\int_{\mathbb{R}^N}f_1(u_n)u_n-\underset{n\to+\infty}{\limsup}\int_{\mathbb{R}^N}f_2(u_n)u_n\\
&=\underset{n\to+\infty}{\liminf}\left[\int_{\mathbb{R}^N}f_1(u_n)u_n-\int_{\mathbb{R}^N}f_2(u_n)u_n\right]\\
&=\underset{n\to+\infty}{\lim}\int_{\mathbb{R}^N}f(u_n)u_n=\left(a+bA^2\right)A^2.
\end{aligned}
\end{equation*}

Now, no matter $A>0$ or $A=0$, with additional help of \eqref{equ2.2} and the following inequality
\begin{equation*}
\int_{\mathbb{R}^N}|\nabla u|^2\leq\underset{n\to+\infty}{\liminf}\int_{\mathbb{R}^N}|\nabla u_n|^2=A^2,
\end{equation*}
we can always conclude that
\begin{equation}\label{equ2.4}
\int_{\mathbb{R}^N}|\nabla u|^2=A^2=\underset{n\to+\infty}{\lim}\int_{\mathbb{R}^N}|\nabla u_n|^2,
\end{equation}
and then
\begin{equation*}
\int_{\mathbb{R}^N}f_2(u)u=\underset{n\to+\infty}{\lim}\int_{\mathbb{R}^N}f_2(u_n)u_n.
\end{equation*}
Noting that $f_2(t)t=2\nu t^2+q(t)$ for all $t\in\mathbb{R}$, where $q(\cdot)$ is a nonnegative continuous function on $\mathbb{R}$ and $\nu>0$. By Fatou's lemma, we have
\begin{equation*}
\underset{n\to+\infty}{\liminf}\int_{\mathbb{R}^N}q(u_n)\geq \int_{\mathbb{R}^N}q(u)~~~~\text{and}~~~~\underset{n\to+\infty}{\liminf}\int_{\mathbb{R}^N}u^2_n\geq \int_{\mathbb{R}^N}u^2,
\end{equation*}
and then
\begin{equation*}
\begin{aligned}
2\nu \int_{\mathbb{R}^N}u^2
&=\int_{\mathbb{R}^N}f_2(u)u-\int_{\mathbb{R}^N}q(u)\geq\underset{n\to+\infty}{\lim}\int_{\mathbb{R}^N} f_2(u_n)u_n-\underset{n\to+\infty}{\liminf}\int_{\mathbb{R}^N}q(u_n)\\
&=2\nu~\underset{n\to+\infty}{\limsup}\int_{\mathbb{R}^N}u^2_n
\geq2\nu~\underset{n\to+\infty}{\liminf}\int_{\mathbb{R}^N}u^2_n\geq2\nu\int_{\mathbb{R}^N}u^2.
\end{aligned}
\end{equation*}
Thus,
\begin{equation}\label{equ2.5}
\underset{n\to+\infty}{\lim}\int_{\mathbb{R}^N}u^2_n=\int_{\mathbb{R}^N}u^2.
\end{equation}
From \eqref{equ2.4}, \eqref{equ2.5} and the fact that $u_n\rightharpoonup u$ in $H^1_r(\mathbb{R}^N)$, we conclude finally that, up to a subsequence, $\{u_n\}^{+\infty}_{n=1}$ converges to $u$ strongly in $H^1_r(\mathbb{R}^N)$ as $n\to+\infty$.~~$\square$

\section{Proof of \autoref{thm 1.1}}\label{section 3}
In this section, we shall deal with the case $a=0$, $b>0$ and $N=2,3$ and give a variational proof of \autoref{thm 1.1}. To be more precise, infinitely many radial solutions are derived in Subsection \ref{subsection 3.1} following the symmetric mountain pass argument and a positive ground state solution is obtained in Subsection \ref{subsection 3.2} by using the mountain pass approach. To avoid distracting the reader from the main line of the proof, some necessary but more technical auxiliary results and their proofs are put into Appendix \hyperref[appendix]{A}.

\subsection{Existence of infinitely many radial solutions}\label{subsection 3.1}

In view of Item $(ii)$ of \autoref{lemma 2.1}, for every $k\in \mathbb{N}^+$, we can define a family of mappings $\Gamma_k$ by
\begin{equation}\label{equ3.1}
  \Gamma_k:=\left\{\gamma\in C(\mathbb{D}_k,H^1_r(\mathbb{R}^N))~|~\gamma~\text{is odd and}~\gamma(\sigma)=\gamma_{0k}(\sigma)~\text{on}~\sigma\in \mathbb{S}^{k-1}\right\}.
\end{equation}
Clearly, $\Gamma_k$ is nonempty since it contains the following mapping
\begin{equation*}
\gamma_k(\sigma):=\left\{
\begin{aligned}
&|\sigma|\gamma_{0k}\left(\frac{\sigma}{|\sigma|}\right),&~~~~&\text{for}~\sigma\in \mathbb{D}_k\setminus\{0\},&\\
&0,&~~~~&\text{for}~\sigma=0.&\\
\end{aligned}
\right.
\end{equation*}
Thus, for every $k\in \mathbb{N}^+$, the symmetric mountain pass value $c_k$ of $\Phi$ defined by
 \begin{equation}\label{equ3.2}
c_k:=\underset{\gamma\in \Gamma_k}{\inf}\underset{\sigma\in \mathbb{D}_k}{\max}~\Phi(\gamma(\sigma))
 \end{equation}
is meaningful. From \eqref{equA.1} and \autoref{thm A.1}, we further conclude the following result.
\begin{lemma}\label{lemma 3.1}
For every $k\in\mathbb{N}^+$, $c_k\geq d_k\geq \rho_0>0$; and thus, $c_k\to+\infty$ as $k\to+\infty$.
\end{lemma}

Now, to obtain infinitely many radial solutions, it is sufficient to prove that, at least for a infinite subset of $\mathbb{N}^+$, the minimax value $c_k$ defined above is a critical value of $\Phi$. For every $k\in\mathbb{N}^+$, by using Ekeland's principle, we can find a Palais-Smale sequence $\{u_n\}^{+\infty}_{n=1}$ at the level $c_k$, that is, a sequence $\{u_n\}^{+\infty}_{n=1}$ that satisfies, as $n\to+\infty$,
\begin{equation}\label{equ3.3}
\Phi(u_n)\to c_k~~~~\text{and}~~~~~\Phi'(u_n)\to 0~\text{in}~(H^1_r(\mathbb{R}^N))^*.
\end{equation}
Observing also \autoref{lemma 2.5}, a first thought one might have is to try to show that $\{u_n\}^{+\infty}_{n=1}$ is bounded in $H^1_r(\mathbb{R}^N)$. However, the nonlinearity $f$ we consider here is so general that this seems difficult to be proved merely under the condition \eqref{equ3.3}. Thus, instead of the above ``general" one, we next try to find a new particular sequence that can turn out to be bounded in $H^1_r(\mathbb{R}^N)$ and also satisfy \eqref{equ3.3}.

Based on the key idea due to Jeanjean \cite{Je97} and further developed by Hirata-Ikoma-Tanaka \cite{Hi10}, we introduce the auxiliary functional
\begin{equation*}
\Psi(\theta,u):=\frac{1}{4}b e^{2(N-2)\theta}\left(\int_{\mathbb{R}^N}|\nabla u|^2\right)^2-e^{N\theta}\int_{\mathbb{R}^N}F(u)
\end{equation*}
on the augmented space $\mathbb{R}\times H^1_r(\mathbb{R}^N)$. Then, for every $k\in \mathbb{N}^+$, we define the minimax value $\overline{c}_k$ of $\Psi$ as follow:
\begin{equation*}
\overline{c}_k:=\underset{\overline{\gamma}\in \overline{\Gamma}_k}{\inf}\underset{\sigma\in \mathbb{D}_k}{\max}\Psi\left(\overline{\gamma}\left(\sigma\right)\right),
\end{equation*}
where
\begin{equation*}
\overline{\Gamma}_k:=\left\{
\begin{aligned}
\overline{\gamma}\in C(\mathbb{D}_k,\mathbb{R}\times H^1_r(\mathbb{R}^N))\left|
\begin{aligned}
&\overline{\gamma}(\sigma)=(\theta(\sigma),\eta(\sigma))~
\text{satisfies}~\\
&(\theta(-\sigma),\eta(-\sigma))=(\theta(\sigma),-\eta(\sigma)), \forall\sigma\in \mathbb{D}_k,\\
&(\theta(\sigma),\eta(\sigma))=(0,\gamma_{0k}(\sigma)), \forall\sigma\in \mathbb{S}^{k-1}.
\end{aligned}
\right.
\end{aligned}
\right\}
\end{equation*}
and $\gamma_{0k}$ is given in Item $(ii)$ of \autoref{lemma 2.1}. It is clear that $\{(0,\gamma)|\gamma\in\Gamma_k\}\subset \overline{\Gamma}_k$ and, arguing as the proofs of Lemma 4.1 and Proposition 4.2 in \cite{Hi10}, we have the following result.

\begin{lemma}\label{lemma 3.3}
For every $k\in \mathbb{N}^+$, $\overline{c}_k$ is well-defined and $\overline{c}_k=c_k$. In addition, there exists a sequence $\{(\theta_n,w_n)\}^{+\infty}_{n=1}\subset \mathbb{R}\times H^1_r(\mathbb{R}^N)$ such that, as $n\to+\infty$,
\begin{itemize}
\item[~~$(i)$] $\theta_n\to0$,
\item[~$(ii)$] $\Psi(\theta_n,w_n)\to c_k$,
\item[$(iii)$] $\partial_u\Psi(\theta_n,w_n)\to0$ in $(H^1_r(\mathbb{R}^N))^*$,
\item[$(iv)$] $\partial_\theta\Psi(\theta_n,w_n)\to 0$.
\end{itemize}
\end{lemma}

As we can see below, by taking full advantages of \autoref{lemma 3.3}, we can prove that the sequence $\{w_n\}^{+\infty}_{n=1}$ is actually bounded in $H^1_r(\mathbb{R}^N)$ and further satisfies \eqref{equ3.3}.

\begin{lemma}\label{lemma 3.4}
Let $\{(\theta_n,w_n)\}^{+\infty}_{n=1}$ be the sequence given by \autoref{lemma 3.3}. Then $\{w_n\}^{+\infty}_{n=1}$ is a bounded Palais-Smale sequence of $\Phi$ at the symmetric mountain pass value $c_k$.
\end{lemma}
\proof For the sake of completeness and clarity, we provide the detailed proof here and divide it into three claims. We shall first prove in Claim 1 that $\{\|\nabla w_n\|^2_2\}^{+\infty}_{n=1}$ is bounded, then conclude the boundedness of $\{\|w_n\|^2_2\}^{+\infty}_{n=1}$ in Claim 2, and finally complete the proof by showing in Claim 3 that the bounded sequence $\{w_n\}^{+\infty}_{n=1}$ is actually a (bounded) Palais-Smale sequence of $\Phi$ at the level $c_k$.

\medskip
\noindent
\textbf{Claim 1.}~{\it The sequence $\{\|\nabla w_n\|^2_2\}^{+\infty}_{n=1}$ is bounded.}

In view of Items $(ii)$ and $(iv)$ of \autoref{lemma 3.3}, we have
\begin{equation*}
\frac{4-N}{4N}b e^{2(N-2)\theta_n}\left(\int_{\mathbb{R}^N}|\nabla w_n|^2\right)^2
=\Psi(\theta_n,w_n)-\frac{1}{N}\partial_\theta\Psi(\theta_n,w_n)
=c_k+o_n(1).
\end{equation*}
Thus, the boundedness of $\{\|\nabla w_n\|^2_2\}^{+\infty}_{n=1}$ follows directly from Item $(i)$ of \autoref{lemma 3.3} and the fact that $N=2,3$.

\medskip
\noindent
\textbf{Claim 2.}~{\it The sequence $\{\|w_n\|^2_2\}^{+\infty}_{n=1}$ is bounded. Then, by Claim 1, $\{w_n\}^{+\infty}_{n=1}$ is bounded in $H^1_r(\mathbb{R}^N)$.}

Arguing by contradiction, let us assume that, up to a subsequence, $\|w_n\|_2\to+\infty$ as $n\to+\infty$. For every $n\in\mathbb{N}^+$, set $t_n:=\|w_n\|^{-2/N}_2$ and $v_n(\cdot):=w_n(t^{-1}_n\cdot)$. Then,
\begin{equation}\label{equ3.4}
t_n\to 0~~~~\text{as}~ n\to+\infty.
\end{equation}
Since $\|v_n\|^2_2=1$ and $\|\nabla v_n \|^2_2=t^{N-2}_n\|\nabla w_n\|^2_2$, we know from Claim 1 and \eqref{equ3.4} that $\{v_n\}^{+\infty}_{n=1}$ is bounded in $H^1_r(\mathbb{R}^N)$. Up to a subsequence, we may assume that $v_n\rightharpoonup v_0$ in $H^1_r(\mathbb{R}^N)$ and $v_n(x)\to v_0(x)$ almost everywhere in $\mathbb{R}^N$ for some $v_0\in H^1_r(\mathbb{R}^N)$. Set $\varepsilon_n:=\|\partial_u\Phi(\theta_n,w_n)\|_{(H^1_r(\mathbb{R}^N))^*}$, then we have
\begin{equation*}
\begin{split}
2\nu e^{N\theta_n}\leq&~ b t^2_n e^{2(N-2)\theta_n}\int_{\mathbb{R}^N}|\nabla w_n|^2\int_{\mathbb{R}^N}|\nabla v_n|^2+2\nu e^{N\theta_n}\int_{\mathbb{R}^N}v^2_n\\
=&~\partial_u\Psi(\theta_n,w_n)[t^N_n w_n]+e^{N\theta_n}\int_{\mathbb{R}^N}(f(v_n)+2\nu v_n)v_n\\
\leq&~ \varepsilon_n t^{N/2}_n(t^{N}_n\|\nabla w_n\|^2_2+1)^{\frac{1}{2}}+e^{N\theta_n}\int_{\mathbb{R}^N}f_1(v_n)v_n,
\end{split}
\end{equation*}
where $f_1$ is the continuous function introduced in the proof of \autoref{lemma 2.5}. Arguing as the proof of \eqref{equ2.3}, we also have that
\begin{equation*}
\underset{n\to+\infty}{\lim}\int_{\mathbb{R}^N}f_1(v_n)v_n=\int_{\mathbb{R}^N}f_1(v_0)v_0.
\end{equation*}
Thus, in view of \eqref{equ3.4}, Claim 1 and Items $(i)$ and $(iii)$ of \autoref{lemma 3.3}, the following inequality holds
\begin{equation*}
0<2\nu\leq\int_{\mathbb{R}^N}f_1(v_0)v_0,
\end{equation*}
which implies that $v_0\not\equiv0$.

On the other hand, let $z\in C^\infty_0(\mathbb{R}^N)\cap H^1_r(\mathbb{R}^N)$ and set $\varphi_n(\cdot):=z(t_n\cdot)$ for every $n\in\mathbb{N}^+$. Obviously, arguing as the proof of \eqref{equ2.1}, we have
\begin{equation*}
\underset{n\to+\infty}{\lim}\int_{\mathbb{R}^N}f(v_n)z=\int_{\mathbb{R}^N}f(v_0)z.
\end{equation*}
Then, by using \eqref{equ3.4}, Claim 1 and Items $(i)$ and $(iii)$ of \autoref{lemma 3.3} again, we obtain that
\begin{equation*}
\begin{aligned}
\left|\int_{\mathbb{R}^N}f(v_0)z\right|=&~\left|e^{N\theta_n}\int_{\mathbb{R}^N}f(v_n)z\right|+o_n(1)\\
\leq&~\left|\partial_u\Psi(\theta_n,w_n)[t^N_n\varphi_n]\right|+ be^{2(N-2)\theta_n}t^2_n\int_{\mathbb{R}^N}|\nabla w_n|^2\left|\int_{\mathbb{R}^N}\nabla v_n\nabla z\right|+o_n(1)\\
\leq&~ \varepsilon_n t^{N/2}_n(t^2_n\|\nabla z\|^2_2+\|z\|^2_2)^{\frac{1}{2}}+Ct^2_n+o_n(1)\to 0.
\end{aligned}
\end{equation*}
Thus,
\begin{equation*}
\int_{\mathbb{R}^N}f(v_0)z=0~~~~\text{for any}~ z\in C^\infty_0(\mathbb{R}^N)\cap H^1_r(\mathbb{R}^N),
\end{equation*}
which implies that $f(v_0)\equiv0$. However, it follows from $(f2)$ that $0$ is an isolated zero point of the continuous function $f$. Thus, in association with Item $(i)$ of \autoref{lemma 2.3}, we have that $v_0\equiv0$, which is a contradiction.

\medskip
\noindent
\textbf{Claim 3.}~{\it The bounded sequence $\{w_n\}^{+\infty}_{n=1}$ is actually a (bounded) Palais-Smale sequence of $\Phi$ at the level $c_k$.}

Obviously, by using Item $(i)$ of \autoref{lemma 3.3} and Claim 2, Items $(ii)$ and $(iii)$ of \autoref{lemma 3.3}, respectively, gives that, as $n\to+\infty$,
\begin{equation*}
\Phi(w_n)\to c_k~~~~\text{and}~~~~\Phi'(w_n)\to 0~\text{in}~(H^1_r(\mathbb{R}^N))^*.
\end{equation*}
Thus, the proof of this lemma is finished.~~$\square$

Now, with the help of \autoref{lemma 3.4}, \autoref{lemma 2.5} and \autoref{lemma 3.1}, we can draw the final conclusion of this subsection.

\medskip
\noindent
\textbf{Conclusion of Subsection \ref{subsection 3.1}.}~~For every minimax value $c_k>0$, let $\{(\theta_n,w_n)\}^{+\infty}_{n=1}$ be the corresponding sequence given by \autoref{lemma 3.3}. In view of \autoref{lemma 3.4} and \autoref{lemma 2.5}, we may assume that, up to a subsequence, $w_n\to u_k$ in $H^1_r(\mathbb{R}^N)$ for some $u_k\in H^1_r(\mathbb{R}^N)$, and then
\begin{equation*}
\Phi(u_k)=c_k~~~~\text{and}~~~~\Phi'(u_k)=0.
\end{equation*}
Thus, for every $k\in\mathbb{N}^+$, the minimax value $c_k$ defined by \eqref{equ3.2} is indeed a critical value of $\Phi$. Now we know from \autoref{lemma 3.1} that, in the case $a=0$ and $N=2,3$, Problem \eqref{problem P} has infinitely many distinct radial solutions for any $b>0$, the energies of which are convergent to positive infinity. This completes the proof of the multiplicity result claimed in \autoref{thm 1.1}.~~$\square$

\subsection{Existence of a positive ground state solution}\label{subsection 3.2}

We first recall that a nontrivial solution $u$ of Problem \eqref{problem P} is said to be a ground state (or least energy) solution if and only if $\Phi(u)=m$, where $m$ is the least energy level defined as follow
\begin{equation*}
m:=\inf_{w\in\mathcal{S}}\Phi(w)~~~~\text{and}~~~~
\mathcal{S}:=\left\{w\in H^1(\mathbb{R}^N)\setminus\{0\}~|~\Phi'(w)=0\right\}.
\end{equation*}

In our case here, that is $a=0$, $b>0$ and $N=2,3$, from the above subsection and Poho\u{z}aev identity, we know that $m$ is well-defined with $0\leq m\leq c_1$. To show the existence of a positive ground state solution, we introduce the following two mountain pass minimax values:
\begin{equation*}
c_{mp}:=\underset{\gamma\in\Gamma}{\inf}~\underset{t\in[0,1]}{\max}\Phi(\gamma(t))~~~~\text{and}~~~~c_{mp,r}:=\underset{\gamma\in\Gamma_r}{\inf}~\underset{t\in[0,1]}{\max}\Phi(\gamma(t)),
\end{equation*}
where
\begin{equation*}
\Gamma:=\{\gamma\in C([0,1],H^1(\mathbb{R}^N))~|~\gamma(0)=0,~\Phi(\gamma(1))<0\}
\end{equation*}
and
\begin{equation*}
\Gamma_r:=\{\gamma\in C([0,1],H^1_r(\mathbb{R}^N))~|~\gamma(0)=0,~\Phi(\gamma(1))<0\}.
\end{equation*}
In view of \autoref{lemma 2.1}, $c_{mp}$ and $c_{mp,r}$ are well-defined satisfying
\begin{equation*}
 0< \rho_0\leq c_{mp}\leq c_{mp,r}<+\infty.
\end{equation*}
Moreover, we have the following result which provides us with some new characterizations of $m$ and plays a important role in obtaining a positive ground state solution.
\begin{lemma}\label{lemma 3.5}
Let $c_1$, $c_{mp,r}$, $c_{mp}$ and $m$ be the values defined as above. Then,
$$c_1=c_{mp,r}=c_{mp}=m.$$
\end{lemma}
\proof  We observe by \autoref{thm B.1} that $c_1$ actually does not depend on the choice of $\gamma_{01}(1)\in H^1_r(\mathbb{R}^N)$ with $\Phi(\gamma_{01}(1))<0$. This fact implies that $c_1= c_{mp,r}$. As a direct consequence of \autoref{thm C.1}, we also know that $c_{mp}\leq m$. Thus,
\begin{equation*}
c_{mp}\leq m\leq c_1= c_{mp,r}.
\end{equation*}
To complete the proof of this lemma, it is sufficient to show that
\begin{equation}\label{equ3.14}
c_{mp,r}\leq c_{mp}+\delta~~~~\text{for any given}~ \delta>0.
\end{equation}

For this purpose, we adopt the argument explored in the proof of Proposition 2.13 of \cite{Fi14}. Let $\rho\in C^\infty_0(\mathbb{R}^N)$ such that $\rho(x)\geq0$ for all $x\in\mathbb{R}^N$ and $\int_{\mathbb{R}^N}\rho(x)dx=1$. For any given $\gamma\in \Gamma$, $\varepsilon>0$ and $t\in[0,1]$, we set
\begin{equation*}
\gamma_\varepsilon(t)(x):=\frac{1}{\varepsilon^N}\int_{\mathbb{R}^N}\rho\left(\frac{x-y}{\varepsilon}\right)\gamma(t)(y)dy,
\end{equation*}
and denote by $\gamma^*_\varepsilon(t)$ the symmetric decreasing rearrangement of $\gamma_\varepsilon(t)$. We will show that $\gamma^*_\varepsilon\in\Gamma_r$ for sufficiently small $\varepsilon>0$ and that \eqref{equ3.14} indeed holds.

 It is clear that the following statements hold.
\begin{description}
  \item[$~~(i)$] For any $\varepsilon>0$ and $t\in[0,1]$, the function $\gamma_\varepsilon(t)$ is of class $C^{\infty}(\mathbb{R}^N)\cap H^1(\mathbb{R}^N)$.
  \item[$~(ii)$] For any $\varepsilon>0$, the mapping $\gamma_\varepsilon:[0,1]\to H^1(\mathbb{R}^N)$ is continuous.
  \item[$(iii)$] $\max_{t\in[0,1]}\|\gamma_\varepsilon(t)-\gamma(t)\|_{H^1}\to0$ as $\varepsilon\to0^+$.
\end{description}
In addition,
\begin{equation*}
\int_{\mathbb{R}^N}|\nabla \gamma^*_\varepsilon(t)|^2\leq\int_{\mathbb{R}^N}|\nabla \gamma_\varepsilon(t)|^2
~~~~\text{and}~~~~\int_{\mathbb{R}^N}F(\gamma^*_\varepsilon(t))=\int_{\mathbb{R}^N}F(\gamma_\varepsilon(t)),
\end{equation*}
which imply that
\begin{equation}\label{equ3.8}
\Phi(\gamma^*_\varepsilon(t))\leq\Phi(\gamma_\varepsilon(t)).
\end{equation}
In view of Item $(iii)$, we have that $\Phi(\gamma_\varepsilon(1))<0$ for sufficiently small $\varepsilon>0$ and
\begin{equation}\label{equ3.7}
\underset{t\in[0,1]}{\max} \Phi(\gamma_\varepsilon(t))\to\underset{t\in[0,1]}{\max}\Phi(\gamma(t))~~~~\text{as}~\varepsilon\to0^+.
\end{equation}
Thus, it follows from \eqref{equ3.8} that
\begin{equation}\label{equ3.9}
\Phi(\gamma^*_\varepsilon(1))<0~~~~\text{for sufficiently small}~\varepsilon>0.
\end{equation}
On the other hand, since $\gamma_\varepsilon(t)\in C^\infty(\mathbb{R}^N)$, we observe by Theorem 5.2 of \cite{Al89} that $\gamma_\varepsilon(t)$ is co-area regular. (For the definition of co-area regularity, we refer readers to Definition 1.2.6 of \cite{Al89}.) In association with Item $(ii)$, we conclude by Theorem 1.4 of \cite{Al89} that
\begin{equation}\label{equ3.10}
\gamma^*_\varepsilon\in C([0,1],H^1_r(\mathbb{R}^N))~~~~\text{for any given}~\varepsilon>0.
\end{equation}
Now, from \eqref{equ3.10}, the fact that $\gamma^*_\varepsilon(0)=0$ and \eqref{equ3.9}, it follows that $\gamma^*_\varepsilon\in\Gamma_r$ for sufficiently small $\varepsilon>0$. Thus, for any given $\delta>0$ and sufficiently small $\varepsilon>0$, by using the definition of $c_{mp,r}$, \eqref{equ3.8} and \eqref{equ3.7}, we have that
\begin{equation*}
c_{mp,r}\leq \underset{t\in[0,1]}{\max}\Phi(\gamma^*_\varepsilon(t))
\leq\underset{t\in[0,1]}{\max}\Phi(\gamma_\varepsilon(t))\leq\underset{t\in[0,1]}{\max}\Phi(\gamma(t))+\delta.
\end{equation*}
Since $\gamma\in \Gamma$ is arbitrary, we know that \eqref{equ3.14} holds which completes the proof.~~$\square$

The reader should be aware that the radial solution of Problem \eqref{problem P} corresponding to $c_1$ given by the above subsection may not be positive (or negative). Thus, even though we have now $c_1=m$ by \autoref{lemma 3.5}, it is still not sufficient to conclude that Problem \eqref{problem P} has a positive ground state solution.

Fortunately, as we can see below, following the argument of the proof of Theorem 6.3 in \cite{Hi10}, we are able to show that Problem \eqref{problem P} has a positive radial solution corresponding to the mountain pass value $c_{mp,r}$. In view of \autoref{lemma 3.5}, this solution turns out to be a positive ground state solution of Problem \eqref{problem P}.

\begin{lemma}\label{lemma 3.6}
Problem \eqref{problem P} has a positive radial solution corresponding to the level $c_{mp,r}$.
\end{lemma}
\proof We use again the auxiliary functional $\Psi(\theta,u)$ introduced in above subsection and define here a new minimax value $\overline{c}_{mp,r}$ of $\Psi$ as follow:
\begin{equation*}
\overline{c}_{mp,r}:=\underset{\overline{\gamma}\in \overline{\Gamma}_r}{\inf}\underset{t\in [0,1]}{\max}\Psi(\overline{\gamma}(t)),
\end{equation*}
where
\begin{equation*}
\overline{\Gamma}_r:=\left\{
\begin{aligned}
\overline{\gamma}\in C([0,1],\mathbb{R}\times H^1_r(\mathbb{R}^N))\left|
\begin{aligned}
&\overline{\gamma}(t)=\left(\theta(t),\eta(t)\right)~
\text{satisfies}\\
&\theta(0)=0=\theta(1),\\
&\eta(0)=0~\text{and}~\Psi(\overline{\gamma}(1))<0.\\
\end{aligned}
\right.
\end{aligned}
\right\}.
\end{equation*}
Obviously, $\{(0,\gamma)|\gamma\in\Gamma_r\}\subset \overline{\Gamma}_r$ and, arguing as the proof of Lemma 4.1 in \cite{Hi10}, we have that $\overline{c}_{mp,r}=c_{mp,r}$.

We observe that $\Phi(u)=\Phi(|u|)$ for all $u\in H^1_r(\mathbb{R}^N)$ and the mapping
\begin{equation*}
|\cdot|:H^1_r(\mathbb{R}^N)\to H^1_r(\mathbb{R}^N),~~~~u\mapsto |u|
\end{equation*}
is continuous with respect to the standard $H^1_r\text{-norm}$. Thus, for any given $n\in\mathbb{N}^+$, without loss of generality, we may assume that there exists a path $\gamma_n\in\Gamma_r$ such that
\begin{equation}\label{equ3.11}
\underset{t\in[0,1]}{\max}\Phi(\gamma_n(t))\leq c_{mp,r}+\frac{1}{n}
\end{equation}
and
\begin{equation}\label{equ3.12}
\gamma_n(t)(x)\geq0~~~~\text{for all}~t\in[0,1]~\text{and}~x\in\mathbb{R}^N.
\end{equation}

Now, based on \eqref{equ3.11} and the fact that $(0,\gamma_n)\in\overline{\Gamma}_r$ for all $n\in\mathbb{N}^+$ and $\overline{c}_{mp,r}=c_{mp,r}$, arguing as the proof of Proposition 4.2 in \cite{Hi10}, we are able to find a sequence $\{(\theta_n,w_n)\}^{+\infty}_{n=1}\subset \mathbb{R}\times H^1_r(\mathbb{R}^N)$ such that, as $n\to+\infty$,
\begin{itemize}
\item[~~$(i)$] $\theta_n\to0$,
\item[~$(ii)$] $\Psi(\theta_n,w_n)\to c_{mp,r}$,
\item[$(iii)$] $\partial_u\Psi(\theta_n,w_n)\to0$ in $(H^1_r(\mathbb{R}^N))^*$,
\item[$(iv)$] $\partial_\theta\Psi (\theta_n,w_n)\to 0$,
\end{itemize}
and, for any given $n\in\mathbb{N}^+$,
\begin{equation}\label{equ3.13}
\text{dist}_{\mathbb{R}\times H^1_r(\mathbb{R}^N)} \left((\theta_n,w_n),\{0\}\times\gamma_n([0,1])\right)\leq \frac{2}{\sqrt{n}}.
\end{equation}
Clearly, by repeating the argument of the proof of \autoref{lemma 3.4}, we know that $\{w_n\}^{+\infty}_{n=1}$ is a bounded Palais-Smale sequence of $\Phi$ at the mountain pass value $c_{mp,r}$. In view of \autoref{lemma 2.5}, we may assume that, up to a subsequence,
\begin{equation*}
w_n\to u_0~~~~\text{in}~H^1_r(\mathbb{R}^N)
\end{equation*}
for some $u_0\in H^1_r(\mathbb{R}^N)$. Thus, $u_0$ is a nontrivial radial solution of Problem \eqref{problem P} satisfying $\Phi(u_0)=c_{mp,r}$. Moreover, \eqref{equ3.12} and \eqref{equ3.13} give that, as $n\to+\infty$,
\begin{equation*}
\left\|(w_n)^-\right\|_2\leq\text{dist}_{\mathbb{R}\times H^1_r(\mathbb{R}^N)} \left((\theta_n,w_n),\{0\}\times\gamma_n([0,1])\right)\to 0,
\end{equation*}
where $w^-(x):=\max\{0,-w(x)\}$. This implies that $(u_0)^-\equiv0$ and, by the maximal principle, $u_0(x)>0$ for all $x\in\mathbb{R}^N$.~~$\square$

\bigskip
\noindent
\textbf{Conclusion of Subsection \ref{subsection 3.2}.}~~In view of \hyperref[lemma 3.5]{Lemmas \ref{lemma 3.5}} and \ref{lemma 3.6}, we know that Problem \eqref{problem P} has a positive ground state solution which is radially symmetric. This completes the proof of the remaining part of \autoref{thm 1.1}.~~$\square$

\section{Proof of \autoref{thm 1.2}}\label{section 4}

In this section, we shall deal with the case $a\geq0,b>0$ and $N\geq5$ and prove \autoref{thm 1.2} by using Clark theorem. Before starting the proof, we give some notations and state the Clark theorem.

Let $X$ be a real Banach space. A subset $B$ of $X$ is said to be symmetric if $u\in B$ implies $-u\in B$. Denote by $\it\Sigma$ the family of closed symmetric subsets of $X$ which do not contain $0\in X$. For any given $B\in\it\Sigma$, the genus $\mathcal{G}(B)$ of $B$ is by definition the smallest integer $k\in\mathbb{N}^+$ for which there exists an odd continuous mapping $\phi:B\to \mathbb{R}^k\setminus\{0\}$, $\mathcal{G}(B)=+\infty$ if no such mapping exists, and $\mathcal{G}(\emptyset)=0$. Now, we state the Clark theorem of the following form.

\begin{theorem}\label{thm 4.1} (\cite{Cl72,He87}) Let $X$ be a real Banach space with norm $\|\cdot\|$ and $I\in C^1(X,\mathbb{R})$. Assume that $I$ is even and bounded from below, satisfies the Palais-Smale condition and $I(0)=0$. Then the following statements hold.
\begin{description}
  \item[$~(i)$] If, for some $k\in\mathbb{N}^+$, there exists $B_k\in \Sigma$ such that
                   \begin{equation*}
                          \mathcal{G}(B_k)\geq k~~~~\text{and}~~~~\sup\nolimits_{u\in B_k}I(u)<0,
                   \end{equation*}
                then $I$ possesses at least $k$ distinct critical points with negative energies.
  \item[$(ii)$] If the assumption of Item $(i)$ holds for every $k\in\mathbb{N}^+$, then $I$ has a sequence of critical points $\{u_k\}^{+\infty}_{k=1}$ satisfying $I(u_k)<0$ for all $k\in\mathbb{N}^+$ and $I(u_k)\to0$ as $k\to+\infty$.
\end{description}
\end{theorem}

\begin{remark}\label{remark 4.1}
If additionally $0\in X$ is the unique critical point of $I$ with zero energy, by virtue of the Palais-Smale condition of $I$, we can conclude further that the sequence of critical points $\{u_k\}^{+\infty}_{k=1}$ given by Item $(ii)$ of \autoref{thm 4.1} converges strongly to zero in $X$. Otherwise, the sequence of critical points $\{u_k\}^{+\infty}_{k=1}$ does not necessarily converge to zero, see \cite{Ka05} for such an example. We also refer readers to \cite{Ch17,Ka05,Li15} for certain improved versions of the Clark theorem and their various applications.
\end{remark}

In the context here, we set $X=H^1_r(\mathbb{R}^N)$, equip $X$ with the standard norm $\|\cdot\|_{H^1_r}$ and let
\begin{equation*}
I(u):=\Phi(u)=\frac{a}{2}\int_{\mathbb{R}^N}|\nabla{u}|^2+\frac{b}{4}\left(\int_{\mathbb{R}^N}|\nabla{u}|^2\right)^2-\int_{\mathbb{R}^N}F(u),
\end{equation*}
where $a\geq0$, $b>0$ and $N\geq5$. Obviously, $I\in C^1(X,\mathbb{R})$ is even and $I(0)=0$. Moreover, $I$ is bounded from below due to Item $(i)$ of \autoref{lemma 2.2} and satisfies the Palais-Smale condition because of Item $(i)$ of \autoref{lemma 2.2} and \autoref{lemma 2.5}. Now, we are ready to prove \autoref{thm 1.2}.

If $a>0$ fixed, for every given $k\in\mathbb{N}^+$, let
\begin{equation*}
b\in(0,b_k)~~~~\text{and}~~~~B_k:=\overline{\gamma}_{0k}(\mathbb{S}^{k-1}),
\end{equation*}
where $b_k>0$ and $\overline{\gamma}_{0k}\in C(\mathbb{S}^{k-1}, H^1_r(\mathbb{R}^N)\setminus\{0\})$ given by Item $(ii)$ of \autoref{lemma 2.2}. Obviously, for $b\in(0,b_k)$, we have
\begin{equation*}
B_k\in\Sigma,~~~~\mathcal{G}(B_k)\geq k~~\text{and}~~\sup\nolimits_{u\in B_k}I(u)<0.
\end{equation*}
Thus, we conclude from Item $(i)$ of \autoref{thm 4.1} that Problem \eqref{problem P} has at least $k$ distinct radial solutions with negative energies for any $b\in(0,b_k)$.

If $a=0$, for every $k\in\mathbb{N}^+$, let
\begin{equation*}
B_k:=\widetilde{\gamma}_{0k}(\mathbb{S}^{k-1}),
\end{equation*}
where $\widetilde{\gamma}_{0k}\in C(\mathbb{S}^{k-1},H^1_r(\mathbb{R}^N)\setminus\{0\})$ given by Item $(iii)$ of \autoref{lemma 2.2}. It is clear that
\begin{equation*}
B_k\in\Sigma,~~~~\mathcal{G}(B_k)\geq k~~\text{and}~~\sup\nolimits_{u\in B_k}I(u)<0.
\end{equation*}
On the other hand, $0\in H^1_r(\mathbb{R}^N)$ is the unique critical point of $I$ with zero energy. Indeed, in this case, any critical point $u$ of $I$ satisfies the following Poho\u{z}aev identity
\begin{equation*}
\frac{N-2}{2N}b\left({\int_{\mathbb{R}^N}}|\nabla{u}|^2\right)^2=\int_{\mathbb{R}^N}F(u).
\end{equation*}
Since $N\geq5$, then
\begin{equation*}
I(u)=-\frac{N-4}{4N}b\left({\int_{\mathbb{R}^N}}|\nabla{u}|^2\right)^2\leq0,
\end{equation*}
where the equality holds if and only if $u=0$. Thus, in view of Item $(ii)$ of \autoref{thm 4.1} and \autoref{remark 4.1}, we conclude that Problem \eqref{problem P} has infinitely many distinct radial solutions $\{u_k\}^{+\infty}_{k=1}$ for any $b>0$. In addition, $\Phi(u_k)<0$ for all $k\in\mathbb{N}^+$ and $u_k\to0$ in $H^1(\mathbb{R}^N)$ as $k\to+\infty$.

\begin{remark}\label{remark 4.2}
These radial solutions $\{u_k\}^{+\infty}_{k=1}$ are actually of class $C^2$ and thus of class $L^\infty$. Unlike the case of the sublinear problems considered in \cite{Ka05,Li15}, we cannot show here that $\|u_k\|_\infty\to0$ as $k\to+\infty$, even though $\{u_k\}^{+\infty}_{k=1}$ converges strongly to zero in $H^1_r(\mathbb{R}^N)$ and converges almost everywhere to zero in $\mathbb{R}^N$ (since $|u(x)|\leq C_N|x|^{\frac{2-N}{2}}\|\nabla u\|_2$ for $u\in \mathcal{D}^{1,2}_r(\mathbb{R}^N)$ and every $x\neq 0$, see Radial Lemma A.III in \cite{Be83-1}). Actually, we have here that
\begin{equation*}
\|u_k\|_\infty>\zeta_0>0~~\text{for all}~k\in\mathbb{N}^+,
\end{equation*}
where $\zeta_0:=\inf\{t>0~|~F(t)>0\}$. It is clear that $\zeta_0$ is well defined because of $(f4)$ and is positive due to $(f2)$. Moreover, $F(t)\leq0$ for all $t\in [-\zeta_0, \zeta_0]$. If $\|u_k\|_\infty\leq\zeta_0$ for some $k\in\mathbb{N}^+$, by using Poho\u{z}aev identity,  we have
\begin{equation*}
0<\frac{N-2}{2N}b\left({\int_{\mathbb{R}^N}}|\nabla{u_k}|^2\right)^2=\int_{\mathbb{R}^N}F(u_k)\leq0,
\end{equation*}
which gives a contradiction.
\end{remark}

\bigskip
We close with some remarks concerning the remaining degenerate cases $a=0,b>0,N=1$ and $a=0,b>0,N=4$.

When $a=0,b>0$ and $N=1$, no compactness can be regained as that in Subsection \ref{subsection 2.2}. Thus, to deal with this case, different methods are needed. The argument in \cite{Je03} (see also Subsection 2.1 of \cite{Fi14}) maybe a possible one.

When $a=0,b>0$ and $N=4$, for any given critical point $u$ of $\Phi$, we know from Poho\u{z}aev identity that $\Phi(u)=0$. This fact means that, in this case, $0\in\mathbb{R}$ is the only potential level at which a nontrivial critical point of $\Phi$ may exists. On the other hand, it was proved in \cite{Lu17} that Problem \eqref{problem P} has only the trivial solution when $b>0$ is large and has nontrivial solutions if $b>0$ is one of a certain positive sequence $\{b_k\}^{+\infty}_{k=1}$ which converges to zero, see Item $(ii)$ of Theorem 1.3 in \cite{Lu17}. Thus, to deal with this case, it seems that a more natural and proper way is to treat Problem \eqref{problem P} as the following nonlinear nonlocal eigenvalue problem
\begin{equation*}
-\left(\int_{\mathbb{R}^4}|\nabla{u}|^2\right)\Delta{u}= \lambda f(u),~~~~(\lambda,u)\in \mathbb{R}\times H^1(\mathbb{R}^4),
\end{equation*}
rather than to look for critical points of $\Phi$. We will consider this nonlinear nonlocal eigenvalue problem as well as its non-autonomous case (from a variational point of view) in a future work.

\bigskip
\section*{Acknowledgment}
\addcontentsline{toc}{section}{Acknowledgment}
The author would like to thank Professor Zhi-Qiang Wang for his valuable comments and constant encouragement.

\begin{appendix}
\section*{Appendix A}\label{appendix}
\addcontentsline{toc}{section}{Appendix A}
\renewcommand{\thesubsection}{A.\arabic{subsection}}
\setcounter{theorem}{0}
\renewcommand{\thetheorem}{A.\arabic{theorem}}
\setcounter{definition}{0}
\renewcommand{\thedefinition}{A.\arabic{definition}}
\setcounter{lemma}{0}
\renewcommand{\thelemma}{A.\arabic{lemma}}
\setcounter{proposition}{0}
\renewcommand{\theproposition}{A.\arabic{proposition}}
\setcounter{remark}{0}
\renewcommand{\theremark}{A.\arabic{remark}}
\setcounter{example}{0}
\renewcommand{\theexample}{A.\arabic{example}}
\setcounter{corollary}{0}
\renewcommand{\thecorollary}{A.\arabic{corollary}}
\setcounter{equation}{0}
\renewcommand{\theequation}{A.\arabic{equation}}

This appendix involves some technical auxiliary results, which are crucial to conduct the variational proof of \autoref{thm 1.1} in Section \ref{section 3}.

\subsection{An auxiliary problem and its result}\label{appendix A}
To obtain infinitely many radial solutions in Subsection \ref{subsection 3.1}, one of the key points is to show that the sequence of potential critical values $\{c_k\}^{+\infty}_{k=1}$ defined by \eqref{equ3.2} converges to positive infinity. Since we now consider the degenerate Kirchhoff problem, the previous local auxiliary problem presented in \cite{Hi10,Lu16} is actually not valid here any more. Fortunately, as we can see below, we can introduce a new one which is still in the spirit of \cite{Hi10} but turns out to be  sufficient to establish \autoref{lemma 3.1}.

Consider $p_0\in (3,5)$ if $N=3$, $p_0\in (3,+\infty)$ if $N=2$ and set
\begin{equation*}
\begin{split}
&h(t):=\left\{
\begin{aligned}
&\max\{\nu t+f(t),0\},~&\text{for}~t\geq0,\\
&\min\{\nu t+f(t),0\},&\text{for}~t<0,\\
\end{aligned}
\right.
\\
&\overline{h}(t):=\left\{
\begin{aligned}
&t^{p_0}\underset{0<\tau\leq t}{\max}\frac{h(\tau)}{\tau^{p_0}},~~~~~~~&\text{for}~t>0,\\
&0,&\text{for}~t=0,\\
&-|t|^{p_0}\underset{t\leq\tau <0}{\max}\frac{h(|\tau|)}{|\tau|^{p_0}},&\text{for}~t<0,\\
\end{aligned}
\right.\\
&\overline{H}(t):=\int^t_0\overline{h}(\tau)d\tau,
\end{split}
\end{equation*}
where $f$ is a function satisfying $(f1)-(f4)$ and $\nu>0$ is the positive constant defined at the very beginning of Section \ref{section 2}. Then the functions $h,\overline{h}$ and $\overline{H}$ satisfy the properties stated in \autoref{lemma A.1} below, whose proof can be found in \cite{Hi10}.
\begin{lemma}\label{lemma A.1}
Let $h,\overline{h}$ and $\overline{H}$ be the functions defined as above, then the following statements hold.
\begin{itemize}
  \item[$~~(i)$] There exists $\delta>0$ such that $h(t)=\overline{h}(t)=0$ for all $t\in[-\delta,\delta]$.
  \item[$~(ii)$] For all $t\in\mathbb{R}$, we have $\frac{1}{2}\nu t^2+F(t)\leq \overline{H}(t)$.
  \item[$(iii)$] For all $t\in \mathbb{R}$, we have $0\leq (p_0+1)\overline{H}(t)\leq\overline{h}(t)t$.
  \item[$(iv)$] The mapping $t\mapsto \overline{h}(t)-\nu t$ satisfies $(f1)-(f4)$.
\end{itemize}
\end{lemma}

Now, the new auxiliary problem can be constructed as follow:
\begin{equation*}
-\left(b\int_{\mathbb{R}^N}|\nabla u|^2\right)\Delta{u}+\nu u= \overline{h}(u),~~u\in H^1_r(\mathbb{R}^N)\setminus\{0\},
\end{equation*}
where $b>0$, $N=2,3$, $\nu>0$ and $\overline{h}\in C(\mathbb{R},\mathbb{R})$ defined as above. Obviously, the new auxiliary problem is nonlocal and the corresponding functional given by
\begin{equation*}
J(u):=\frac{1}{4}b\left(\int_{\mathbb{R}^N}|\nabla u|^2\right)^2+\frac{1}{2}\nu\int_{\mathbb{R}^N}u^2-\int_{\mathbb{R}^N}\overline{H}(u)
\end{equation*}
is of class $C^1(H^1_r(\mathbb{R}^N))$. Moreover, as stated in the next lemma, the functional $J$ has the symmetric
mountain pass geometry and satisfies the Palais-Smale compactness condition.

\begin{lemma}\label{lemma A.2}
The functional $J$ satisfies \autoref{lemma 2.1} and the Palais-Smale compactness condition.
\end{lemma}
\proof We know from Item $(ii)$ of \autoref{lemma A.1} that, when $a=0$, $b>0$ and $N=2,3$,
\begin{equation}\label{equA.1}
\Phi(u)\geq J(u)~~~~\text{for all}~u\in H^1_r(\mathbb{R}^N).
\end{equation}
Thus, the odd continuous mapping $\gamma_{0k}$ given by Item $(ii)$ of \autoref{lemma 2.1} is also valid here. In addition, in view of Item $(iv)$ of \autoref{lemma A.1}, it is clear that the functional $J$ satisfies Item $(i)$ of \autoref{lemma 2.1} by rechoosing $r_0>0$ and $\rho_0>0$ smaller.

Thanks to Item $(iii)$ of \autoref{lemma A.1} and the fact that $p_0>3$, we can show in a standard way that every Palais-Smale sequence of $J$ is bounded in $H^1_r(\mathbb{R}^N)$. Thus, in view of Item $(iv)$ of \autoref{lemma A.1}, the Palais-Smale compactness condition of $J$ follows directly from  \autoref{lemma 2.5}.~~$\square$

Now, we define the symmetric mountain pass values of the functional $J$ as follow:
 \begin{equation*}
d_k:=\underset{\gamma\in \Gamma_k}{\inf}\underset{\sigma\in \mathbb{D}_k}{\max}~J(\gamma(\sigma)),
 \end{equation*}
where $\Gamma_k$ is given by \eqref{equ3.1} and $k\in\mathbb{N}^+$, and conclude the main result of this subsection.
\begin{theorem}\label{thm A.1}
For every $k\in \mathbb{N}^+$, the level $d_k$ is actually a critical value of $J$ and $d_k\geq\rho_0>0$. In addition, $d_k\to +\infty$ as $k\to+\infty$.
\end{theorem}
\proof Since we have \autoref{lemma A.2}, this result can be proved in the exact same way of Lemma 3.2 of \cite{Hi10}. We omit the details here.~~$\square$

\subsection{Arc-wise connectedness}\label{appendix B}

The main aim of this subsection is to show the arc-wise connectedness of the set $O_r$ defined below under suitable assumptions, which is necessary to the proof of \autoref{lemma 3.5} and seems also to be interesting by itself.
\begin{theorem}\label{thm B.1}
Assume that $a\geq0,b>0$ fixed, $N=2,3$ and $f$ satisfies $(f1)-(f4)$. Let $O_r:=\{u\in H^1_r(\mathbb{R}^N)~|~\Phi(u)<0\}$. Then $O_r$ is arc-wise connected, that is, for any given $u_1,u_2\in O_r$, there exists a continuous path $\gamma$ in $O_r$ joining $u_1$ and $u_2$.
\end{theorem}

To prove \autoref{thm B.1}, motivated partly by the proof of Lemma 6.1 in \cite{Hi10}, for $R\geq1$, $s\in[0,1]$ and $t\geq1$, we set
\begin{equation*}
\eta(R,s,t;x):=\pi(R,s;t^{-1}R^{-2}|x|),
\end{equation*}
where
\begin{equation*}
\pi(R,s;r):=
\left\{
\begin{aligned}
&0,&&\text{for}~r\in[0,R],\\
&\zeta(r-R),&&\text{for}~r\in[R,R+1)],\\
&\zeta,&&\text{for}~r\in[R+1,R+1+sR],\\
&\zeta(R+2+sR-r),&&\text{for}~r\in[R+1+sR,R+2+sR],\\
&0,&&\text{for}~r\in[R+2+sR,+\infty),
\end{aligned}
\right.
\end{equation*}
and $\zeta>0$ is given in $(f4)$ satisfying $F(\zeta)>0$. We will see that $\eta(\underline{R},1,T;x)\in O_r$ for sufficiently large $\underline{R}$ and $T$, and there exist continuous curves joining $u_i~(i=1,2)$ and $\eta(\underline{R},1,T;x)$ in $O_r$. Clearly, this proves our \autoref{thm B.1}.

We start with the following result.
\begin{lemma}\label{lemma B.1}
There exist $C_0,C_1,C_2>0$ (which are independent of $R$, $s$, $t$ and $\theta$) and $R^*\geq1$ such that
\begin{description}
  \item[$~~(i)$] $\Phi(\eta(R,s,1;x))\leq C_0 R^{3N-1}$ for all $R\geq 1$ and $s\in[0,1]$.
  \item[$~(ii)$] $\Phi(\theta\eta(R,0,1;x))\leq C_1R^{3N-1}$ for all $R\geq1$ and $\theta\in[0,1]$.
  \item[$(iii)$] $\Phi(\eta(R,1,t;x))\leq -C_2t^N R^{3N}$ for all $R\geq R^*$ and $t\geq1$.
\end{description}
\end{lemma}
\proof For reader's convenience, we provide the detailed proof here. For $R\geq1$, $s\in[0,1]$ and $t\geq1$, by some direct computations, we have
\begin{equation*}
\begin{split}
&\int_{\mathbb{R}^N}|\nabla\eta(R,s,t;x)|^2dx\\
=&~\omega_{N-1}(tR^2)^{N-2}\left(\int^{R+1}_R+\int^{R+2+sR}_{R+1+sR}\right)|\pi_r(R,s;r)|^2r^{N-1}dr\\
=&~\frac{1}{N}\omega_{N-1}\zeta^2(tR^2)^{N-2}\left[(R+1)^N-R^N+(R+2+sR)^N-(R+1+sR)^N\right]
\end{split}
\end{equation*}
and
\begin{equation*}
\begin{split}
&\int_{\mathbb{R}^N}F(\eta(R,s,t;x))dx\\
=&~\omega_{N-1}(tR^2)^N\left(\int^{R+1}_R+\int^{R+1+sR}_{R+1}+\int^{R+2+sR}_{R+1+sR}\right)F(\pi(R,s;r))r^{N-1}dr\\
\geq&~\frac{1}{N}\omega_{N-1}F(\zeta)(tR^2)^N\left[(R+1+sR)^N-(R+1)^N\right]\\
&~~~~~~~~~~-\frac{1}{N}\omega_{N-1}M(tR^2)^N\left[(R+1)^N-R^N+(R+2+sR)^N-(R+1+sR)^N\right],
\end{split}
\end{equation*}
where $\omega_{N-1}$ is the surface area of the unit sphere in $\mathbb{R}^N$ and $M:=\max_{l\in[0,\zeta]}|F(l)|>0$. Obviously, for $R\geq1$ and $s\in[0,1]$, the following inequalities hold
\begin{equation*}
\begin{split}
(R+1)^N-R^N&\leq (2^N-1)R^{N-1},\\
(R+2+sR)^N-(R+1+sR)^N&\leq3^{N-1}(2^N-1)R^{N-1},\\
(R+1+sR)^N-(R+1)^N&\geq s^NR^N.
\end{split}
\end{equation*}
Thus, there exist positive constants $C_3,C_4$ and $C_5$ independent of $R$, $s$ and $t$ such that
\begin{equation}\label{equB.1}
\int_{\mathbb{R}^N}|\nabla\eta(R,s,t;x)|^2dx \leq C_3 t^{N-2}R^{3N-5}
\end{equation}
and
\begin{equation*}
  \int_{\mathbb{R}^N}F(\eta(R,s,t;x))dx\geq t^N(C_4s^N R^{3N}-C_5 R^{3N-1}).
\end{equation*}
Then a positive constant $C_0>0$ can be found, which is also independent of $R$, $s$ and $t$, such that
\begin{equation}\label{equB.2}
\begin{split}
\Phi(\eta(R,s,t;x))&\leq \frac{1}{2}aC_3t^{N-2}R^{3N-5}+\frac{1}{4}b(C_3t^{N-2}R^{3N-5})^2\\
&~~~~~~~~~~~~~~~~~~~~~~~~~~~~~~-t^N(C_4 s^NR^{3N}-C_5 R^{3N-1})\\
&\leq t^N(C_0 R^{3N-1}-C_4 s^N R^{3N}).
\end{split}
\end{equation}

Now, we are ready to complete the proof this lemma. Obviously, by letting $t=1$ in \eqref{equB.2}, Item $(i)$ holds. Arguing as the proof of \eqref{equB.2}, a positive constant $C_1>0$ exists, which is independent of $R\geq1$ and $\theta \in [0,1]$, such that
\begin{equation*}
\begin{split}
\Phi(\theta\eta(R,0,1;x))&\leq \frac{1}{2}aC_3R^{3N-5}+\frac{1}{4}b(C_3R^{3N-5})^2+\frac{1}{N}\omega_{N-1}MR^{2N}\left[(R+2)^N-R^N\right]\\
&\leq C_1 R^{3N-1}.
\end{split}
\end{equation*}
This is Item $(ii)$. Finally, by using \eqref{equB.2} again, there exist $C_2>0$, which is independent of $R$ and $t$, and $R^*\geq1$ such that
\begin{equation*}
\Phi(\eta(R,1,t;x))\leq -C_2t^N R^{3N}~~~~\text{for all}~R\geq R^*~\text{and}~t\geq1.
\end{equation*}
Thus, Item $(iii)$ holds.~~$\square$

We know from Item $(iii)$ of \autoref{lemma B.1} that, under the assumptions of \autoref{thm B.1}, $\eta(\underline{R},1,T;x)\in O_r$ for sufficiently large $\underline{R},T$. For any given $u_i\in O_r~(i=1,2)$, we next try to join $u_i$ and $\eta(\underline{R},1,T;x)$ in $O_r$. Without loss of generality, we may assume that $u_i$ has compact support and $\text{supp}~u_i(x)\subset B(0,L^3_i)$ for some constant $L_i>0$. Now, we consider the following curves:
\begin{equation*}
\begin{aligned}
&\gamma_{i,1}: [L_i, \underline{R}]\to H^1(\mathbb{R}^N);&R&\mapsto u_i((L_i/R)^3x),\\
&\gamma_{i,2}: [0, 1]\to H^1(\mathbb{R}^N);&\theta&\mapsto u_i((L_i/\underline{R})^3x)+\theta\eta(\underline{R},0,1;x),\\
&\gamma_{i,3}: [0, 1]\to H^1(\mathbb{R}^N);&s&\mapsto u_i((L_i/\underline{R})^3x)+\eta(\underline{R},s,1;x),\\
&\gamma_{i,4}: [1, T]\to H^1(\mathbb{R}^N);&t&\mapsto u_i((L_i/\underline{R})^3x)+\eta(\underline{R},1,t;x),\\
&\gamma_{i,5}: [0, 1]\to H^1(\mathbb{R}^N);&\theta&\mapsto (1-\theta)u_i((L_i/\underline{R})^3x)+\eta(\underline{R},1,T;x).
\end{aligned}
\end{equation*}
Since $u_i,\eta(\underline{R},1,T;x)\in H^1_r(\mathbb{R}^N)$, by joining these curves, we get a continuous path in $H^1_r(\mathbb{R}^N)$ joining $u_i$ and $\eta(\underline{R},1,T;x)$. Actually, we will further see that, with suitable choices of $\underline{R}$ and $T$, our path defined above is indeed included in $O_r$.

\begin{lemma}\label{lemma B.2}
For any given $u_i\in O_r$ $(i=1,2)$, the following statements hold.
\begin{description}
  \item[$~~(i)$] $\Phi(u_i((L_i/R)^3x))\leq\Phi(u_i)(R/L_i)^{3N}\leq \Phi(u_i)<0$ for all $R\in[L_i,+\infty)$.
  \item[$~(ii)$] There exists $R_i\geq1$ such that, for all $s\in[0,1]$, $t\geq1$, $\theta\in[0,1]$ and $R\geq R_i$,
  \begin{eqnarray*}
    &\Phi(u_i((L_i/R)^3x)+\theta\eta(R,0,1;x))< 0, \\
    &\Phi(u_i((L_i/R)^3x)+\eta(R,s,1;x))< 0,\\
    &\Phi(u_i((L_i/R)^3x)+\eta(R,1,t;x))<0.
  \end{eqnarray*}
  \item[$(iii)$] There exists $T_i\geq 1$ such that, for all $R\geq R_i$, $\theta\in[0,1]$ and $t\geq T_i$,
    \begin{equation*}
    \Phi((1-\theta)u_i((L_i/R)^3x)+\eta(R,1,t;x))< 0.
  \end{equation*}
\end{description}
\end{lemma}
\proof $(i)$ Since $u_i\in O_r$ and $N=2,3$, we can easily obtain that, for all $R\geq L_i$,
\begin{equation*}
\begin{split}
&~~~~~\Phi(u_i((L_i/R)^3x))\\
&= \frac{1}{2}a(R/L_i)^{3(N-2)}\int_{\mathbb{R}^N}|\nabla u_i|^2+\frac{1}{4}b(R/L_i)^{6(N-2)}\left(\int_{\mathbb{R}^N}|\nabla u_i|^2\right)^2-(R/L_i)^{3N}\int_{\mathbb{R}^N}F(u_i)\\
&\leq\Phi(u_i)(R/L_i)^{3N}\leq \Phi(u_i)<0.
\end{split}
\end{equation*}
Thus, Item $(i)$ holds.

\smallskip
$(ii)$ We note that, for all $R\geq\max\{1,L_i\}$, $s\in[0,1]$ and $t\geq1$,
\begin{equation}\label{equB.3}
\text{supp}~u_i((L_i/R)^3x)\cap\text{supp}~\eta(R,s,t;x)=\emptyset.
\end{equation}
Thus, with additional help of Item $(i)$ of this lemma, Item $(ii)$ of \autoref{lemma B.1} and \eqref{equB.1}, we have that, for all $R\geq\max\{1,L_i\}$ and $\theta\in[0,1]$,
\begin{equation}\label{equB.4}
\begin{split}
&~\Phi(u_i((L_i/R)^3x)+\theta\eta(R,0,1;x))\\
=&~\Phi(u_i((L_i/R)^3x))+\Phi(\theta\eta(R,0,1;x))\\
&~\hspace{20mm}+\frac{1}{2}b(R/L_i)^{3(N-2)}\int_{\mathbb{R}^N}|\nabla u_i|^2\cdot\theta^2\int_{\mathbb{R}^N}|\nabla\eta(R,0,1;x)|^2\\
\leq&~ \Phi(u_i)(R/L_i)^{3N}+ C_1R^{3N-1}+\left(\frac{1}{2}bC_3L^{3N-5}_i\int_{\mathbb{R}^N}|\nabla u_i|^2\right) (R/L_i)^{6N-11}.
\end{split}
\end{equation}
Similarly, we conclude from Item $(i)$ of \autoref{lemma B.1} that, for all $R\geq\max\{1,L_i\}$ and $s\in[0,1]$,
\begin{equation}\label{equB.5}
\begin{split}
&~\Phi(u_i((L_i/R)^3x)+\eta(R,s,1;x))\\
\leq&~\Phi(u_i)(R/L_i)^{3N}+ C_0R^{3N-1}+\left(\frac{1}{2}bC_3L^{3N-5}_i\int_{\mathbb{R}^N}|\nabla u_i|^2\right) (R/L_i)^{6N-11}.
\end{split}
\end{equation}
Also, by using Item $(iii)$ of \autoref{lemma B.1}, we obtain that, for all $R\geq\max\{R^*,L_i\}$ and $t\geq1$,
\begin{equation}\label{equB.6}
\begin{split}
&~\Phi(u_i((L_i/R)^3x)+\eta(R,1,t;x))\\
\leq&~ \Phi(u_i)(R/L_i)^{3N}-\left[C_2 R^{3N}- \left(\frac{1}{2}bC_3L^{3N-5}_i\int_{\mathbb{R}^N}|\nabla u_i|^2\right) (R/L_i)^{6N-11}\right]t^N.
\end{split}
\end{equation}
Since $\Phi(u_i)<0$ and $N=2,3$, we now conclude from \eqref{equB.4}, \eqref{equB.5} and \eqref{equB.6} that a large enough positive constant $R_i$ exists such that Item $(ii)$ holds.

\smallskip
$(iii)$ Obviously, there exists $M_i>0$ independent of $\theta$ such that, for all $\theta\in[0,1]$,
\begin{equation*}
\int_{\mathbb{R}^N}F((1-\theta)u_i)\geq -M_i.
\end{equation*}
Then, for all $R\geq L_i$ and $\theta\in[0,1]$, the following inequality holds
\begin{equation*}
\begin{split}
&~\Phi((1-\theta)u_i((L_i/R)^3x))\\
\leq&~ \left[\frac{1}{2}a\int_{\mathbb{R}^N}|\nabla u_i|^2+\frac{1}{4}b\left(\int_{\mathbb{R}^N}|\nabla u_i|^2\right)^2+M_i\right](R/L_i)^{3N}=:\overline{M}_i(R/L_i)^{3N}.
\end{split}
\end{equation*}
In view of \eqref{equB.3}, Item $(iii)$ of \autoref{lemma B.1} and \eqref{equB.1}, we have that, for all $R\geq R_i$, $\theta\in[0,1]$ and $t\geq1$,
\begin{equation*}
\begin{split}
&~\Phi((1-\theta)u_i((L_i/R)^3x)+\eta(R,1,t;x))\\
\leq &~\left[\overline{M}_i+\left(\frac{1}{2}bC_3L^{3(N-2)}_i\int_{\mathbb{R}^N}|\nabla u_i|^2\right)t^{N-2}- C_2 L^{3N}_i t^N \right](R/L_i)^{3N}.
\end{split}
\end{equation*}
Thus, a large enough positive constant $T_i$ exists such that Item $(iii)$ holds.~~$\square$

\bigskip
\noindent
\textbf{Proof of \autoref{thm B.1}.}~~Let
\begin{equation*}
\underline{R}:=\max\{R_1,R_2\}~~~~\text{and}~~~~T:=\max\{T_1,T_2\},
\end{equation*}
where $R_i$ and $T_i$ $(i=1,2)$ are given by \autoref{lemma B.2}. Obviously, by \autoref{lemma B.2}, we have
\begin{equation*}
\gamma_{i,1}([L_i, \underline{R}]),\gamma_{i,2}([0, 1]),\gamma_{i,3}([0, 1]),\gamma_{i,4}([1, T]),\gamma_{i,5}([0, 1])\subset O_r.
\end{equation*}
Thus, $u_i$ $(i=1,2)$ and $\eta(\underline{R},1,T;x)$ are connected by a continuous path in $O_r$.  Since the positive constants $\underline{R}$ and $T$ are chosen uniformly in $u_i~(i=1,2)$, the proof of \autoref{thm B.1} is completed.~~$\square$

\begin{remark}\label{remark B.2}
We would like to point out that, under the assumptions of \autoref{thm B.1}, the new set $O:=\{u\in H^1(\mathbb{R}^N)~|~\Phi(u)<0\}$ is also arc-wise connected. Indeed, for any given $u_1,u_2\in O$, one can define a continuous path in $H^1(\mathbb{R}^N)$ exactly as above, which is not radially symmetric in general (since $u_1$ or $u_2$ may not be radially symmetric) but joins $u_1$ and $u_2$ in $O$.
\end{remark}

\subsection{Existence of optimal path}\label{appendix C}

We shall prove in this subsection the following result, which claims the existence of an optimal path for any given nontrivial solution of Problem \eqref{problem P} in the case $a=0,b>0$ and $N=2,3$, and is also needed for the proof of \autoref{lemma 3.5}.

\begin{theorem}\label{thm C.1}
Assume that $a=0,b>0$ fixed, $N=2,3$ and $f$ satisfies $(f1)-(f4)$. Then, for any given nontrivial solution $w\in H^1(\mathbb{R}^N)$ of Problem \eqref{problem P}, there exists an optimal path $\gamma$, that is, $\gamma\in\Gamma$ satisfying
\begin{equation*}
w\in \gamma([0,1])~~~~\text{and}~~~~\underset{t\in [0,1]}{\max}\Phi(\gamma(t))=\Phi(w).
\end{equation*}
\end{theorem}
\proof For any given nontrivial solution $w\in H^1(\mathbb{R}^N)$ of Problem \eqref{problem P}, without loss of generality, we will find a curve $\gamma:[0,L]\to H^1(\mathbb{R}^N)$ such that
\begin{equation*}
\gamma(0)=0,~~\Phi(\gamma(L))<0,~~w\in\gamma([0,L])~~\text{and}~~\underset{t\in [0,L]}{\max}\Phi(\gamma(t))=\Phi(w).
\end{equation*}

When $N=3$, the construction of such a curve is rather simple. Actually, we can easily verify that the curve defined by
\begin{equation*}
\gamma(t)(x):=\left\{
\begin{aligned}
&~0,&t=0,~~~~\\
&w(x/t),&t\in(0,L],
\end{aligned}
\right.
\end{equation*}
is the desired one by choosing $L>0$ sufficiently large.

When $N=2$, the situation becomes more complicated. As we can see below, with suitable choices of $t_0\in(0,1)$, $t_1\in (1,+\infty)$ and $\theta_1>1$, the curve $\gamma$ constituted of the three pieces defined below gives a desired one:
\begin{equation*}
\begin{aligned}
&\gamma_1(\theta)(x):=
\left\{
\begin{aligned}
&~~0,\\
&\theta w(t^{-1}_0\theta^{-2}x),
\end{aligned}
\right.
&&\begin{aligned}
&\theta=0,\\
&\theta\in(0,1],
\end{aligned}\\
&\gamma_2(t)(x):=w(x/t),& &t\in[t_0,t_1],\\
&\gamma_3(\theta)(x):=\theta w(t^{-1}_1\theta^{-2} x), &&\theta\in[1,\theta_1].
\end{aligned}
\end{equation*}
Such a construction is motivated by an idea of Jeanjean and Tanaka, which was developed in \cite{Je02} to deal with Problem \eqref{problem Q}. But the curve constructed here is slightly different from the one presented in \cite{Je02}, since we are now faced with some difficulties which are cased by the vanishing of the term $\int_{\mathbb{R}^2}|\nabla{u}|^2$ and the loss of a limit $\lim_{t\to0}f(t)/t$.

For any given $t\in(0,+\infty)$, let
\begin{equation*}
g^t(\theta):=\frac{1}{4}b\theta^4\|\nabla w\|^4_2-\theta^4t^2\int_{\mathbb{R}^2}F(\theta w),~~~~\theta\in[0,+\infty).
\end{equation*}
Obviously, $g^t(\theta)=\Phi(\theta w(t^{-1}\theta^{-2}x))$ when $\theta\in(0,+\infty)$ and
\begin{equation*}
\frac{d}{d\theta}g^t(\theta)
=\theta^3\left[b\|\nabla w\|^4_2-t^2\left(4\int_{\mathbb{R}^2}F(\theta w)+\int_{\mathbb{R}^2}f(\theta w)\theta w\right)\right].
\end{equation*}
We choose $t_0\in(0,1)$ sufficiently small such that
\begin{equation}\label{equC.1}
b\|\nabla w\|^4_2-t^2_0\left(4\int_{\mathbb{R}^2}F(\theta w)+\int_{\mathbb{R}^2}f(\theta w)\theta w\right)>0~~~~\text{for all}~\theta\in[0,1].
\end{equation}
Since $w\in H^1(\mathbb{R}^2)$ is a nontrivial solution of Problem \eqref{problem P}, we have
\begin{equation*}
\int_{\mathbb{R}^2}F(w)=0~~~~\text{and}~~~~\int_{\mathbb{R}^2}f(w)w=b\|\nabla w\|^4_2>0.
\end{equation*}
Thus, a suitable constant $\theta_1>1$ exists independently of $t_0$ such that
\begin{equation*}
4\int_{\mathbb{R}^2}F(\theta w)+\int_{\mathbb{R}^2}f(\theta w)\theta w>0~~~~\text{for all}~\theta\in[1,\theta_1],
\end{equation*}
For this constant $\theta_1$, we choose $t_1>1$ sufficiently large such that, for all $\theta\in[1,\theta_1]$,
\begin{equation}\label{equC.2}
b\|\nabla{w}\|^4_2-t^2_1\left(4\int_{\mathbb{R}^2}F(\theta w)+\int_{\mathbb{R}^2}f(\theta w)\theta w\right)\leq -\frac{2}{\theta^4_1-1}b\|\nabla{w}\|^4_2.
\end{equation}
Now, we know from \eqref{equC.1} that $\Phi(\gamma_1(\cdot))$ is increasing in $\theta\in[0,1]$ and takes its maximal at $\theta=1$. Since $\int_{\mathbb{R}^2}F(w)=0$, we have $\Phi(\gamma_2(\cdot))=\Phi(w)=\frac{1}{4}b\|\nabla{w}\|^4_2$ for all $t\in[t_0,t_1]$. Finally, in view of \eqref{equC.2}, $\Phi(\gamma_3(\cdot))$ is decreasing in $\theta\in[1,\theta_1]$ and
\begin{equation*}
\begin{split}
\Phi(\theta_1w(t^{-1}_1\theta^{-2}_1x))=g^{t_1}(\theta_1)&=g^{t_1}(1)+\int^{\theta_1}_1\frac{d}{d\theta}g^{t_1}(\theta)d\theta\\
&\leq \frac{1}{4}b\|\nabla{w}\|^4_2-\int^{\theta_1}_1\frac{2\theta^3}{\theta^4_1-1}b\|\nabla{w}\|^4_2 d\theta\\
&=-\frac{1}{4}b\|\nabla{w}\|^4_2<0.
\end{split}
\end{equation*}
Therefore, we get the desired curve for the case $N=2$.~~$\square$
\begin{remark}\label{remark C.1}
From the construction process above, we know that, for a positive (radial, respectively) solution of Problem \eqref{problem P}, the corresponding optimal path $\gamma$ actually can be chosen such that $\gamma(t)(x)>0$ for all $x\in\mathbb{R}^N$ and $t\in(0,1]$ ($\gamma(t)$ is radially symmetric for all $t\in(0,1]$, respectively).
\end{remark}
\end{appendix}


{\small

\begin{thebibliography}{100}
\bibitem{Ad00} S. Adachi,  K. Tanaka,
{ Trudinger type inequalities in $\mathbb{R}^N$ and their best exponents},
{\it Proc. Amer. Math. Soc.} {128} (7) (2000), 2051--2057.
\bibitem{Al89} F.J. Almgren, E.H. Lieb,
{ Symmetric decreasing rearangement is sometimes continuous},
J. Amer. Math. Soc. {2}(4) (1989) 683--773.
\bibitem{Az12} A. Azzollini,
{ The elliptic Kirchhoff equation in ${\mathbb{R}^N}$ perturbed by a local nonlinearity},
Differential Integral Equations  {25} (5-6) (2012) 543--554.
\bibitem{Az15} A. Azzollini,
{ A note on the elliptic Kirchhoff equation in ${\mathbb{R}^N}$ perturbed by a local nonlinearity},
Commun. Contemp. Math. {17} (2015) 1450039. 5 pages.
\bibitem{Az11} A. Azzollini, P. d'Avenia, A. Pomponio,
{ Multiple critical points for a class of nonlinear functions},
Ann. Mat. Pura Appl. {190} (2011) 507--523.
\bibitem{Be83-3} H. Berestycki, T. Gallouet, O. Kavian,
{ Equations de champs scalaires euclidens non lineaires dans le plan},
 C. R. Acad. Sci. Paris, Serie I Math {297} (1983) 307--310.
\bibitem{Be83-1} H. Berestycki, P.L. Lions,
{ Nonlinear scalar field equations I: Existence of a ground state},
Arch. Rat. Mech. Anal. {82} (1983) 313--346.
\bibitem{Be83-2} H. Berestycki, P.L. Lions,
{ Nonlinear scalar field equations II: Existence of infinitely many solutions},
Arch. Rat. Mech. Anal. {82} (1983) 347--375.
\bibitem{Ch17} S.W. Chen, Z.L. Liu, Z.-Q. Wang,
{A variant of Clark's theorem and its applications for nonsmooth functionals without the Palais-Smale condition},
SIAM J. Math. Anal. {49} (2017) 446--470.
\bibitem{Cl72} D.C. Clark,
{ A variant of the Lusternik-Schnirelman theory},
{\it Indiana Univ. Math. J.} 22 (1972-1973) 65--74.
\bibitem{Fi14} G.M. Figueiredo, N. Ikoma and J.R.S. J\'{u}nior,
{ Existence and concentration result for the Kirchhoff type equations with general nonlinearities},
{\it Arch. Ration. Mech. Anal.} {\textbf{213}} (2014), 931--979.
\bibitem{He87} H.P. Heinz,
{ Free Ljusternik-Schnirelman theory and the bifurcation diagrams of certain singular nonlinear systems},
{\it J. Differential Equations}  66 (1987) 263--300.
\bibitem{Hi10} J. Hirata, N. Ikoma and K. Tanaka,
{ Nonlinear scalar field equations in $\mathbb{R}^N$: mountain pass and symmetric mountain pass approaches},
{\it Topol. Methods Nonlinear Anal.} {35} (2010) 253--276.
\bibitem{Je97} L. Jeanjean,
{ Existence of solutions with prescribed norm for semilinear elliptic equations},
{\it Nonlinear Anal.} {28} (1997), 1633--1659.
\bibitem{Je02} L. Jeanjean,  K. Tanaka,
{ A remark on least energy solutions in $\mathbb{R}^N$},
{\it Proc. Amer. Math. Soc.} {131} (8) (2002), 2399--2408.
\bibitem{Je03} L. Jeanjean, K. Tanaka,
{ A note on a mountain pass characterization of least energy solutions},
Adv. Nonlinear Stud. {3} (2003) 445--455.
\bibitem{Ka05} R. Kajikiya,
{ A critical point theorem related to the symmetric mountain pass lemma and its applications to elliptic equations},
{\it J. Funct. Anal.} 225 (2005) 352--370.
\bibitem{Li15} Z.L. Liu, Z.-Q. Wang,
{ On Clark's theorem and its applications to partially sublinear problems},
Ann. Inst. H. Poincar\'{e} Anal. Non-Lin\'{e}aire 32 (2015), 1015--1037.
\bibitem{Lu15} S.-S. Lu,
{ Signed and sign-changing solutions for a Kirchhoff-type equation in bounded domains},
J. Math. Anal. Appl. 432 (2015), 965--982.
\bibitem{Lu17} S.-S. Lu,
{ An autonomous Kirchhoff-type equation with general nonlinearity in $\mathbb{R}^N$}, Nonlinear Anal. RWA 34 (2017) 233--249.
\bibitem{Lu16} S.-S. Lu,
{ Multiple solutions for a Kirchhoff-type equation with general nonlinearity}, Adv. Nonlinear Anal. (2016), http://dx.doi.org/10.1515/anona-2016-0093.
\bibitem{St77} W.A. Strauss,
{ Existence of solitary waves in higher dimensions}, Comm. Math. Phys. 55 (1977), 149--162.
\bibitem{Og90} T. Ogawa,
{A proof of Trudinger's inequality and its application to nonlinear Schr\"{o}dinger equations}, Nonlinear Anal.  14 (9) (1990), 765--769.
\end{thebibliography}
}
\end{document}